\documentclass[11pt,a4paper]{article}
\usepackage{amssymb,amsmath,amsthm}
\usepackage{epsfig}

\def\Q{\mathbb{Q}}
\def\N{\mathbb{N}}
\def\C{\mathbb{C}}
\def\Z{\mathbb{Z}}
\def\D{\mathbb{D}}

\def\0{ {\mbox{\tt 0}} }
\def\1{ {\mbox{\tt 1}} }

\newcommand \iter[2]{#1^{\circ #2}}
\def\upper{\mathcal{A}}
\def\orb{\text{orb}}

\newtheorem{theorem}{Theorem}[section]
\newtheorem{proposition}[theorem]{Proposition}
\newtheorem{lemma}[theorem]{Lemma}
\newtheorem{definition}[theorem]{Definition}
\newtheorem{corollary}[theorem]{Corollary}
\newtheorem{deflem}[theorem]{Definition and Lemma}

\newenvironment{thm}[3]{\noindent{\bf #1 #2 (#3)} \itshape}

\def\lineclear
     {\rule{0pt}{0pt}\nopagebreak\par\nopagebreak\noindent}

\title{On the Structure of Abstract Hubbard Trees and the Space of Abstract Kneading Sequences of Degree Two}
\author{Alexandra Kaffl}
\date{}

\begin{document}

\maketitle

\begin{abstract}

 One of the fundamental properties of the Mandelbrot set is that the set of postcritically finite
parameters is structured like a tree. We extend this result to the set of quadratic kneading
sequences and show that this space contains no irrational de\-co\-rations. Along the
way, we prove a combinatorial analogue to the correspondence principle of dynamic and parameter
rays. Our key tool is to work simultaneously with the two equivalent combinatorial concepts of
Hubbard trees and kneading sequences.

\end{abstract}

\section{Introduction}
For studying the structure of the Mandelbrot set $\mathcal{M}$ (the connectedness locus of complex
quadratic polynomials), combinatorial model spaces proved to be very helpful. In this context, the most important combinatorial concepts are
\emph{Hubbard trees, kneading sequences} and \emph{external angles}. In the \emph{Orsay
Notes}~\cite{orsay}, Douady and Hubbard showed that for polynomials whose critical points are all preperiodic the
dynamical behavior is completely encoded in the so called \emph{Hubbard tree}. This result was extended to the set of all postcritically finite polynomials by Poirier \cite{poirier}. The Hubbard tree of a postcritically finite polynomial is a topological
tree obtained by connecting all points on the critical orbits within the filled-in Julia set.
Penrose introduced a model for unicritical polynomials, which is based on gluing together sequences
over a two letter alphabet \cite{penrose}. The concept of external angles was used by Thurston to define
\emph{laminations} on the circle. The quotient of the circle by the lamination is homeomorphic to
the Julia set if and only if the Julia set is locally connected. The same is true for the lamination model of
the Mandelbrot set $\mathcal{M}$.

In this paper, we investigate the space of quadratic kneading sequences as they were defined by Bruin and Schleicher in \cite{BS}. A subset of this space is a model for $\mathcal{M}$. More
precisely, we consider the \emph{parameter tree}, which is the set
\[ 
	\Sigma^{\star\star}:=\{ \nu\in \{\0,\1\}^\N:\;\nu_1=\1,\,\nu \mbox{ non-periodic} \} \cup
	\left\{\nu \mbox{ is $\star$-periodic}\right\}
\] 
together with a partial order $<$. (We use the convention that $\N=\{1,2,\ldots\}$.) A sequence is
called \emph{$\star$-periodic} if it is of the form $\overline{1\nu_2\ldots\nu_n\star}$ or equals
$\overline{\star}$. Our main result extends the \emph{Branch Theorem} for postcritically finite parameters in $\mathcal{M}$  \cite[Proposition XXII.3]{orsay} to the space
$(\Sigma^{\star\star},<)$. It positively answers the central open question in \cite{BS}, that is
where in the parameter tree branching can happen.

\smallskip
The Branch Theorem is one of the main results of Douday's and Hubbard's study of \emph{nervures}. Its importance lies in that it describes
the structure of the set of postcritically finite polynomials in $\mathcal{M}$ and that it is a key step in proving that local
connectivity of $\mathcal{M}$ implies density of hyperbolicity.
More precisely, the Branch Theorem says that there is a partial order on the set of postcritically finite parameters of $\mathcal{M}$ such that two parameters  $c',c''$ either can be compared or there is a maximal parameter $c$ such that $c<c'$ and $c<c''$. 

\smallskip
Before we introduce this order and state the Branch Theorem in a way that justifies this name, we recall some facts about the Mandelbrot set.
By the Riemann mapping theorem, there is a conformal isomorphism $\Phi_{\mathcal{M}}:\C\setminus \mathcal{M}\longrightarrow\C\setminus \overline{\D}$. The image of a radial line in $\C\setminus\overline{\D}$ under $\Phi^{-1}_{\mathcal{M}}$ is called an \emph{external ray}. External rays foliate the complement of $\mathcal{M}$.
A \emph{hyperbolic component} of $\mathcal{M}$ is a connected component of $\{c\in\mathcal{M}:z\mapsto z^2+c \mbox{ has an attracting periodic orbit}\}$, a \emph{Misiurewicz point} is  a parameter $c\in\mathcal{M}$ for which the critical point of $z\mapsto z^2+c$ is preperiodic under iteration. Note that for us, preperiodic means strictly preperiodic.
 It is well known that for any hyperbolic component $W$, the \emph{multiplier map} $\mu: \overline{W}\longrightarrow \overline{\D}$ is a homeomorphism and that at every parameter $c\in \partial W$ with $\mu(c)=e^{2\pi i\theta}$, $\theta\in\Q/\Z$, two external rays with periodic angles (periodic under angle-doubling) are landing. The point $c_W\in \partial W$ with $\mu(c_W)=1$ is called the \emph{root} of $W$. The \emph{wake} of $W$ is the open region of $\C$ that is separated from the origin by $c_W$ and the two external rays landing at $c_W$. The $\frac{p}{q}$-\emph{subwake} of $W$ is the open region in $\C$ separated from $0$ by a parameter $c\in\partial W$ with $\mu(c)=e^{2\pi i \frac{p}{q}}$ and the two external rays landing at $c$ ($\frac{p}{q}$ in lowest terms). 
Every Misiurewicz point is the landing point of a finite, positive number of external rays with preperiodic angles. A \emph{subwake} of a Misiurewicz point $c$ is any of the open regions in $\C$ separated from $0$ by $c$ and two adjacent external rays landing at $c$. The \emph{wake} of $c$ is the union of all its subwakes.

Lau and Schleicher introduce in \cite{lausch} an \emph{order} on the set of postcritically finite parameters of $\mathcal{M}$ and define \emph{combinatorial arcs}: for two hyperbolic components or Misiurewicz points $A,B$, set $A\prec B$ if two external rays landing at (the root of) $A$ separate $B$ from the origin.  The combinatorial arc $[A,B]$ is the collection of all hyperbolic components and Misiurewicz points $C$ such that $A\prec C\prec B$, together with $A$ and $B$. With this notation we can state the Branch Theorem as follows.

\smallskip
\noindent{\bf Theorem (Douady, Hubbard)}\lineclear
   {\it Let $A\not=B$ be two hyperbolic components or Misiurewicz points. Then there is
    a unique
    hyperbolic component or Misiurewicz point $C$ such that $[W_0,A]\cap[W_0,B]=[W_0,C]$, where $W_0$ denotes the main cardioid of $\mathcal{M}$. \qed}

  \smallskip
In analogy to this, there is a partial order $<$ on the set $\Sigma^{\star\star}$ (c.f. Definition \ref{DEForder}), that allows us to define for any two pre- or
$\star$-periodic kneading sequences $\nu\leq\nu^\prime$ the \emph{arc} $[\nu,\nu^\prime]$ that connects them:
\[
	[\nu,\nu^\prime]:=\{\mu\in\Sigma^{\star\star}:\, \mu \mbox{ pre- or $\star$-periodic such that }
	\nu\leq\mu\leq\nu^\prime \}.
\]
 With this notation our main result reads as follows:

\begin{theorem}[Branch Theorem]\label{THMbranchThm}\lineclear
    Let $\nu,\tilde{\nu}$ be $\star$-periodic or preperiodic kneading sequences.
    Then there is a unique kneading sequence $\mu$ such that exactly one of the following
    holds:
    \begin{enumerate}\renewcommand{\labelenumi}{(\roman{enumi})}
        \item   $[\overline{\star},\nu]\cap[\overline{\star},\tilde{\nu}]=[\overline{\star},\mu]$,
                        where $\mu$ is either $\star$-periodic or preperiodic.
        \item    $[\overline{\star},\nu]\cap[\overline{\star},\tilde{\nu}]=[\overline{\star},\mu]\setminus\{\mu\}$,
                    where $\mu$ is $\star$-periodic and primitive.
    \end{enumerate}
\end{theorem}

The discussion in \cite{BS} shows that the second case of the Branch Theorem \ref{THMbranchThm}
occurs if and only if at least one of the two sequences $\nu,\tilde{\nu}$ is not generated by a
polynomial. Since every postcritically finite parameter in $\mathcal{M}$ defines a unique element in
$\Sigma^{\star\star}$, our result provides a combinatorial-topological proof of the classical
Branch Theorem. Furthermore, it shows that the overall structure of kneading sequences generated by
postcritically finite quadratic polynomials extends to the fairly larger set of abstract ($\star$-
and preperiodic) kneading sequences. 

Our key tool in proving this theorem is the parallel use of two equivalent combinatorial
concepts: kneading sequences and Hubbard trees. This technique was introduced by Bruin and
Schleicher in \cite{BS}. In this monograph, they define Hubbard trees in an abstract way, namely as topological trees with dynamics that meet certain requirements. Unlike Hubbard trees in the sense of Douady and Hubbard, these Hubbard trees are not imbedded in the complex plane. The advantage of this
abstract definition is that one has to deal with less
information. Moreover, it is more natural when considering the
relation between Hubbard trees and kneading sequences: in~\cite[Chapter 3]{BS}, Bruin and Schleicher show that there is a
bijection between the set of (equivalence classes of) Hubbard trees and $\star$- or preperiodic kneading sequences. The
order defined on $\Sigma^{\star\star}$ is based on comparing the Hubbard trees associated to
the considered kneading sequences.

In \cite{kafflThesis} we extend this abstract notion of Hubbard trees and kneading sequences to the setting of unicritical polynomials of degree $d\geq 2$. In this case, too, there is a bijection between the set of (equivalence classes) of Hubbard trees and the set of $\star$- and preperiodic kneading sequences. We define a partial order on the space $\Sigma^\star_d$ of kneading sequences of degree $d$ in an analogous way and extend the Branch Theorem to  $\Sigma^{\star}_d$.

\smallskip
A central part of \cite{BS} is to investigate the properties of Hubbard trees and to classify
Hubbard trees which are generated by quadratic polynomials. For proving the Branch Theorem, we
continue the investigation of the structural properties of Hubbard trees. In Theorem
\ref{THMexistenceBifuracation}, we show the existence of dynamical bifurcation points in a
Hubbard tree $T$. More precisely, if $T$ contains a characteristic point $x$ of
itinerary $\tau$, then it also contains a characteristic point $y$ whose itinerary is a
bifurcation sequence of $\tau$ such that $]x,y[$ contains no characteristic point. With Lemmas
\ref{LEMcharPointsAreForced} and \ref{LEMcharPointAfterBranchPointNotForced} we extend the
\emph{Orbit Forcing Lemma} \cite[Lemma 6.2]{BS}. These lemmas allow us to compare the
characteristic points in two given Hubbard trees $T$ and $T'$. In particular, we can find the
characteristic points in $T$ and $T'$ which have the same itinerary. The itinerary of the point
with this property closest to the critical value gives rise to the branch point $\mu$ of
Theorem~\ref{THMbranchThm}.

\smallskip Our paper is structured the following way: in Chapter \ref{SECdef} we repeat all relevant
definitions introduced in~\cite{BS} and state their main results. In
Chapter~\ref{SECstructureHT} we investigate the structure of Hubbard trees and in Chapter
\ref{SECparamPlane} we derive important properties of the space of kneading sequences and prove the
Branch Theorem \ref{THMbranchThm}.

\bigskip
{\bf Acknowledgements ---} The author would like to thank Dierk Schlei\-cher, Markus F\"orster,
G\"unter Rottenfu\ss er and Johannes R\"uckert for many fruitful discussions and helpful
suggestions. We would like to express our special thanks to the referee for valuable comments and for pointing out a shorter argument for the proof of Theorem \ref{THMexistenceBifuracation}. We are grateful to the Institut Henri Poincar\'e in Paris for its
hospitality, to the Marie-Curie Fellowship Association for financing our stay at the IHP, and to
International University Bremen and the Konrad-Adenauer-Stiftung for their financial support.


\section{Hubbard Trees and Kneading Sequences}\label{SECdef}

In this paper, we focus on two equivalent combinatorial tools used to describe the
dynamical behavior of quadratic polynomials, namely Hubbard trees and kneading sequences. We use the
terminology introduced in~\cite{BS}. To make our paper more self-contained we state the most
important definitions and results of~\cite{BS} in this section.

\subsection{Hubbard Trees}

\begin{definition}[Tree; Branch-, Inner- and Endpoints]\lineclear
    A \emph{tree} is a connected, compact metric space which can be written as a finite union of
    arcs and does not contain a simple closed curve.

    \smallskip
    Let $T$ be a tree and $x\in T$. If the number of components of $T\setminus\{x\}$ equals one, then
    $x$ is called an \emph{endpoint}, if this number is two, then $x$ is called an \emph{inner point} and if
    it is larger than two, $x$ is a \emph{branch point} of $T$.

    \smallskip
    We denote by $[x,y]\subset T$ the unique closed arc in $T$ connecting $x$ and $y$. The arc without
    its endpoints is denoted by $]x,y[$, and the arc containing $x$ but not $y$ by $[x,y[$.
\end{definition}

\smallskip
\begin{definition}[Hubbard Tree]\label{DEFabstrHubbardTree}\lineclear
    \emph{\cite[Definition 3.2]{BS}} A \emph{Hubbard tree} $(T,f)$ is a topological tree $T$ together with a continuous and surjective map
    \mbox{$f:T\to T$} and with a point $c_0$ such that the following hold:
    \newcounter{prop}
    \begin{list}{(\alph{prop})}{\usecounter{prop}  \setlength{\labelwidth}{3ex}
    \setlength{\itemsep}{1ex} \setlength{\parsep}{0ex}}

        \item $f$ is at most 2-to-1,
        \item $f$ is locally 1-to-1 on $T\setminus\{c_0\}$,
        \item all endpoints of $T$ are in $\orb_f(c_0)$,
        \item $c_0$ is periodic or preperiodic,
        \item if $x\not=y$ are \emph{marked points},
            then $\exists\;n\geq 0$ such that $c_0\in\iter{f}{n}([x,y])$.
     \end{list}
     The point $c_0$ is called the \emph{critical point}; a point $p\in T$ is \emph{marked} if it is either a branch point of $T$ or it is contained in the \emph{critical orbit} $\orb_f(c_0)$. The set of all marked points is denoted by $V$. The \emph{critical value} $c_1$ of a Hubbard tree is the image of the critical point $c_0$, i.e., $c_1:=f(c_0)$.
\end{definition}

Observe that this definition of  Hubbard trees is not
equivalent to the definition of Douady and Hubbard in~\cite{orsay}, but it generalizes theirs. We do not require an embedding of $(T,f)$ into the plane, in particular, we do not specify any cyclic order on the arms of branch points of $T$.

\begin{definition}[Equivalent Trees, Arms and Precritical Points]\lineclear
    We say that two Hubbard trees $(T,f)$ and $(T^\prime,f^\prime)$ are \emph{equivalent} if there is a
    bijection $\varphi$ between the set of marked points $V\subset T$ and $V'\subset T'$ which commutes with the dynamics and which has the property that for all $x,y\in V$, $]x,y[\,\cap V=\emptyset$ if and only if $]\varphi(x),\varphi(y)[\,\cap V'=\emptyset$.
    
    \smallskip
    For any point $z\in T$, we call a connected component of $T\setminus\{z\}$ a \emph{global
    arm} of $z$. We define the \emph{local arm} of a global arm $G$ as the direct limit
    $\displaystyle\lim_{\longrightarrow} G\cap U$, where $U$
    is an open neighborhood of $z$. We regard a local arm as a sufficiently small interval $]z,p[\,\subset G$. The \emph{$i$-th image} $\iter{f}{i}(L)$ of a local arm $L$ of $z$ is the local arm at $\iter{f}{i}(z)$ such that its representative intersects $\iter{f}{i}(]z,p[)$, where $]z,p[$ represents the local arm $L$. A local arm  $L$ is \emph{periodic} if $z$ is $n$-periodic and there is an $j>0$ such that $\iter{f}{jn}(L)=L$. The smallest positive number $j$ with this property is called the \emph{period} of $L$.

    \smallskip
     A \emph{precritical point} $\xi_k$ of step $k$ is a point in $T$
    such that $k>0$ is the smallest number with $\iter{f}{k}(\xi_k)=c_1$; we write {\sc step}$(\xi_k)=k$. By this definition,
    $c_0$ is always a precritical point of step $1$, while $c_1$ is a precritical point if and only if $c_0$ is periodic. In this case, {\sc step}$(c_1)$ equals
    the period of $c_0$.
\end{definition}

\begin{remark}(Basic Properties of Hubbard Trees){\bf.}\label{REMbasicProperties}\lineclear
    It is easy to show that the critical value $c_1$ is an endpoint of the tree and consequently the critical point is never a branch point. Thus, no branch
    point maps onto the critical point and all branch points are either periodic or preperiodic~\cite[Lemma 3.3]{BS}.

    Any arc bounded by two precritical points $\xi_k,\tilde{\xi}_k$ with equal step $k$
    contains a precritical point $\xi$ with {\sc step}$(\xi)<k$ because otherwise the
    homeomorphism $\iter{f}{(k-1)}|_{[\xi_k,\tilde{\xi}_k]}$ would map both points $\xi_k,\,\tilde{\xi}_k$ to $c_0$.

    Whenever we speak of the period of a periodic point $p$ we mean its exact period, that is the
    smallest integer $n>0$ such that $\iter{f}{n}(p)=p$.
\end{remark}

\begin{deflem}[Characteristic Point]\lineclear \emph{\cite[Definition \& Lemma 4.1]{BS}}
    Let $\{z_1,\ldots,z_n=z_0\}$ be a periodic orbit in the Hubbard tree $(T,f)$ disjoint from the critical
    orbit.\footnote{For a Hubbard tree with preperiodic critical point we need in addition that no itinerary $\tau(z_i)$ 
        (defined in~\ref{DEFitineraryKS}) coincides with any itinerary of the endpoints of
        $T$, c.f.\ the Hubbard tree associated to $z\mapsto z^2+i.$}

    The $z_i$ can be relabelled in a unique way such that $f(z_i)=z_{i+1}$ (indices are taken modulo $n$) and
    such that one global arm of $z_1$ contains $c_1$ and one other arm contains
   $c_0$ and $z_2,\ldots,z_n$. The point $z_1$ (after relabelling) is called the
    \emph{characteristic point} of the given periodic orbit. \qed
\end{deflem}

If we say that $z\in T$ is a characteristic point, then we mean that $z$ is periodic and $z$ is the characteristic point of $\orb(z)$. Besides precritical points and points in $\orb_f(c_0)$, the characteristic points are the most important
points in the Hubbard tree. They carry the dynamically relevant information. Note that sometimes we
refer to $c_1$ as the characteristic point of the critical orbit.

\begin{figure}[ht] \setlength{\unitlength}{1cm}\begin{center}
      \begin{picture}(11,5)
      	\put(0,3){\text{characteristic point}}
            \put(0,1.5){\epsfig{file=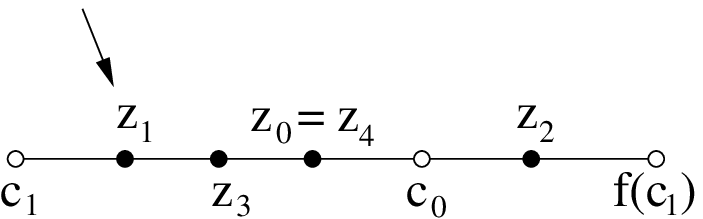, width=4.5cm}}
             \put(5.5,0){\epsfig{file=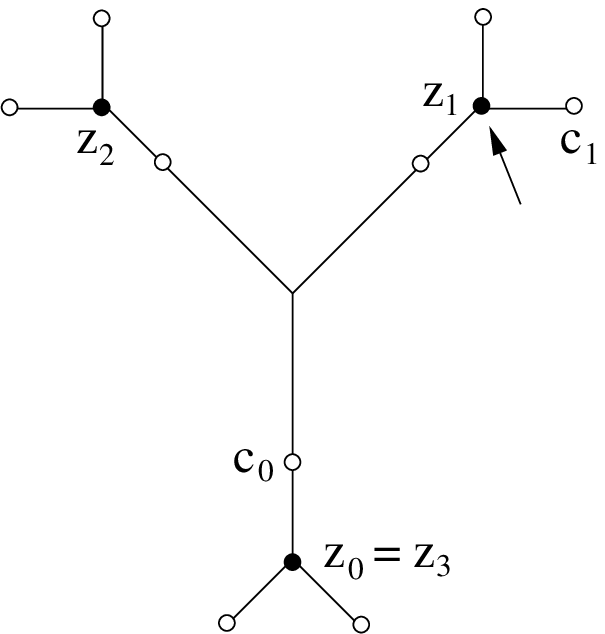, width=4.5cm}}
             \put(8.5,2.7){\text{characteristic}}\put(9.2,2.3){\text{point}}
        \end{picture} \end{center}

        \caption{Characteristic points in the Hubbard tree of the Airplane ($c\approx-1.75$, left) and of the $\frac{1}{3}$-bifurcation of the Rabbit ($c\approx -0.032+0.791i$).}
    \end{figure}

\subsection{Kneading Sequences}

In the second part of this section, we recall some important facts about kneading sequences. We set
$\Sigma^1:=\{\nu=(\nu_i)_{i\geq 1}\in\{\0,\1\}^\N:\,\nu_1=\1  \}$.

\begin{definition}[Kneading Sequences]\label{DEFkneadingSequence}\lineclear
    We say a sequence $\nu\in \{\0,\1,\star\}^\N$ is \emph{$\star$-periodic of period $n$} if it is of the form
    $\nu=\overline{\1\ldots\nu_{n-1}\star}$ with $\nu_i\in\{\0,\1\}$ or equals $\overline{\star}$ for $n=1$.
   A \emph{kneading sequence} is an element of the space
    \[ \Sigma^\star:=\Sigma^1 \cup\,\left\{\nu\in\{\0,\1,\star\}^\N : \,\nu \text{ is
        $\star$-periodic}\right\}.\]
\end{definition}

\smallskip
This notation of kneading sequences differs from the one given in~\cite[Chapter 6]{lausch}, where
only such elements of $\Sigma^\star$ are called kneading sequences which are realized by quadratic
polynomials. Definition \ref{DEFkneadingSequence} also includes such symbolic sequences which are not generated by any quadratic polynomial. We will mainly work with a subset of $\Sigma^{\star}$,
namely with the set
\[\Sigma^{\star\star}:=\Sigma^{\star}\setminus\{\nu \mbox{ is periodic but not $\star$-periodic}\}.\]

\smallskip
Combinatorial objects closely related to kneading sequences are the so called \emph{internal
addresses}. An internal address is a finite or infinite sequence of natural numbers
$1=S_0\rightarrow S_1\rightarrow \cdots \rightarrow S_k \rightarrow \cdots$ such that $S_k<S_{k+1}$
for all $k\in\N$. We call an element $S_k$ an \emph{entry}. 
For any kneading sequence $\nu$ we define the $\rho$-map by
\[
        \rho_\nu: \N \rightarrow \N\cup\{\infty\}, \;  \rho_\nu(n)=
	        \left\{\begin{array}{ll} \inf\{k>n:\, \nu_k\not=\nu_{k-n}\} & \mbox{if existing} \\
	        					\infty & \mbox{otherwise} \end{array}\right..
\]
The \emph{internal address} of $\nu$ is the set $\orb_{\rho_\nu}(1)\setminus\{\infty\}$. In fact, there is a 1-to-1 correspondence between internal addresses and the elements of $\Sigma^1$: each of the two can be calculated from the other via the $\rho$-map  
(c.f.~\cite[Chapter 6]{lausch}). This is not true for $\Sigma^{\star\star}$, but for each element
of $\Sigma^{\star\star}$ there is a unique internal address associated to it.
Suppose that the center of the hyperbolic component $W\subset\mathcal{M}$ generates the kneading sequences $\nu$. Then the entries $S_i$ of the internal address of $\nu$ corresponds to the periods of hyperbolic components $W_i$ which lie on the combinatorial arc $[W_0,W]$ such that $[W_i,W]$ contains no hyperbolic component of period smaller than or equal to the periods of $W$ or $W_i$, see \cite{lausch}. (Here again $W_0$ denotes the main cardioid.)

\smallskip
 Let $\nu$ be any $\star$-periodic kneading sequence of period $n>1$. By replacing the $\star$ by $\0$ everywhere or by $\1$ everywhere,
 we obtain two periodic sequences without the symbol $\star$. Exactly one of those has $n$ as entry in its internal address. This one will be denoted by $\upper(\nu)$ and called the \emph{upper kneading
 sequence} of $\nu$. The other one, $\overline{\upper}(\nu)$, is called the \emph{lower kneading sequence} of $\nu$. 
 
 Suppose that $\nu$ is generated by the Hubbard tree $(T,f)$. Then $\overline{\upper}(\nu)$ corresponds to the limit sequence of the itineraries $\tau(x_n)$ for $n\to\infty$, where $x_n\in T$ such that $x_n\to c_1$ as $n\to\infty$ (itineraries are defined in Definition \ref{DEFitineraryKS}). If $\nu$ comes
 from a polynomial, we can interpret upper and lower kneading sequences also the following way:
there is a hyperbolic component $W$ in $\mathcal{M}$ such that $\nu$ is the kneading sequence generated by the angle $\theta$ of an external ray landing at the root of $W$. The
 lower kneading sequence corresponds to the limit of kneading sequences generated by external angles which converge
 to $\theta$ outside the wake of $W$. The upper kneading sequence corresponds to the limit taken over external angles
 converging to $\theta$ within the wake of $W$.

 We can also go the other way round: given any periodic sequence $\tau$ of period $n>1$ not containing a
 $\star$, we obtain a $\star$-periodic sequence $\upper^{-1}_{ln}(\tau)$
 of period $ln$ by replacing the $j(ln)$-th entry of $\tau$ by the symbol $\star$ for all $j\in\N$.

\smallskip
 Let us further introduce two special kinds of kneading sequences.
 Given any kneading sequence $\nu\in\Sigma^{\star\star}$ with internal address $1\rightarrow S_1\rightarrow \cdots
\rightarrow S_k \rightarrow \cdots$, then $\nu^k$ denotes the $\star$-periodic
sequence corresponding to the internal address $1\rightarrow S_1\rightarrow \cdots \rightarrow
S_k$. Any such sequence is called a \emph{truncated sequence} of $\nu$. Observe that $\nu^k$ is the
unique $\star$-periodic sequence which coincides with $\nu$ for the first $S_k-1$ entries and has
period $S_k$.

For $q>0$, we define the \emph{$q$-th bifurcation sequence $B_q(\mu)$} of a periodic or $\star$-periodic
kneading sequence $\mu=\overline{\mu_1\ldots\mu_n}$ as follows:
    \begin{eqnarray*}
     B_q(\mu):=&\overline{(\mu_1\ldots\mu_n)^{q-1}\mu_1\ldots\mu_{n-1}\mu'_n}\,,
       		 &\text{if } \mu_n\not=\star,\\
       B_q(\mu):=&\overline{(\mu_1\ldots\mu_n')^{q-1}\mu_1\ldots\mu_{n-1}\star}\,,
       & \text{if } \mu_n=\star,
    \end{eqnarray*}
    where in the first case $\mu'_n:=1-\mu_n$, and in the second case $\mu_n'$ is chosen such that $\upper(\mu)=\overline{\mu_1\ldots\mu_n'} $.
 The expression
 $(\mu_1\ldots\mu_n)^{q-1}$ means that the word $\mu_1\ldots\mu_n$ is repeated $(q-1)$-times.

We call a $\star$-periodic kneading sequence $\mu$ \emph{primitive} if $\mu$ and $\overline{\upper}(\mu)$ have the same period. For such sequences, we define the \emph{$q$-th backward bifurcation sequence $\overline{B}_q(\mu)$}  by
\[ 
	\overline{B}_q(\mu):=\overline{(\mu_1\ldots\mu'_n)^{q-1}\mu_1\ldots\mu_{n-1}\star}\,,
	\mbox{ where } \overline{\upper}(\mu)=\overline{\mu_1\ldots\mu'_n}.
\]

If $\nu$ is generated by the hyperbolic component $W$ then the kneading sequence $B_q(\nu)$
 corresponds to all bifurcation components of $W$ at internal angle $\frac{p}{q}$
 (in lowest terms). We will see later that a backward bifurcation corresponds to a bifurcation into a
 non-admissible subtree of $\Sigma^{\star\star}$.

\subsection{Admissibility}\label{subsection_Admissibility}

Recall that the critical point $c_0$ decomposes any Hubbard tree $T$ into at most two components.
Let $T_\1$ be the component containing the critical value and let $T_\0$ be the other component. $T_{\0}$ is empty if and only if $T$ is an $n$-od and all points on the periodic critical orbit are endpoints of $T$. In this case, it follows that for $n>2$, the branch point $b$ is fixed and that for $n=2$, there is a fixed inner point $b\in\,]c_0,c_1[$. In both cases, $b$ equals the $\alpha$-fixed point of $T$. Note that $f$ maps $[b,\iter{f}{i}(c_0)]$ homeomorphically onto $[b,\iter{f}{i+1}(c_0)]$ for all $i\in\N_0$. 

\begin{definition}[Itinerary, Kneading Sequence of $(T,f)$]\label{DEFitineraryKS}\lineclear
    Let $(T,f)$ be a Hubbard tree and $z\in T$. The \emph{itinerary} of $z$ is the infinite sequence
    $\tau(z)=(\tau_i(z))_{i>0}$ defined by
    $$ \tau_i(z)=\left\{\begin{array}{ll} \0 & \text{if } \iter{f}{(i-1)}(z) \in T_{\0}, \\
                                    \star & \text{if } \iter{f}{(i-1)}(z)=c_0, \\
                                    \1 & \text{if } \iter{f}{(i-1)}(z) \in T_{\1}.
                    \end{array}\right.$$
    The \emph{kneading sequence $\nu=(\nu_i)_{i>0}$ of the Hubbard tree $T$} is the itinerary of
    its critical value $c_1$.
\end{definition}

The expansivity condition implies that periodic marked points have the same period as their
itinerary. One fundamental result of~\cite{BS} is the following theorem.
\begin{theorem}[Hubbard Trees vs.\ Kneading Sequences]\label{THMkneadingSequAreRealized}\lineclear \emph{\cite[Theorem
3.11]{BS}}
    Every $ \star$-periodic or preperiodic kneading sequence in $\Sigma^\star$ is realized by a
    Hubbard tree that is unique up to equivalence. \qed
\end{theorem}

Since, by definition, every Hubbard tree generates a unique kneading sequence, we have a 1-to-1
correspondence between the set of pre- and $\star$-periodic kneading sequences and the set of
Hubbard trees (up to equivalence). By the triple $(T,f,\nu)$ we mean the Hubbard tree $(T,f)$
together with its associated kneading sequence $\nu$. The Hubbard tree associated to
$\overline{\star}$ is a point.

\smallskip Observe that there is no 1-to-1 correspondence between Hubbard trees of postcritically finite polynomials
(i.e., Hubbard trees in the sense of ~\cite{orsay}) and the set of $\star$-periodic and preperiodic
kneading sequences. Indeed, because of symmetries in the Mandelbrot set there are several polynomials sharing the same kneading sequence. For example, both the rabbit and the
anti-rabbit have the kneading sequence $\overline{\1\1\star}$. On the other hand, there are
kneading sequences which are not realized at all. An example is the kneading sequence
$\overline{\1\0\1\1\0\star}$, which was discovered by Penrose (amongst others). There are many
kneading sequences which are not generated by polynomials; in fact, we can generate infinitely many such sequences from the kneading sequence associated to any primitive hyperbolic component in $\mathcal{M}$, see~\cite{BS}. We discuss this fact in the beginning of Chapter
\ref{SECparamPlane} in more detail.

\smallskip
One of the key steps in determining which sequences are generated by quadratic polynomials is to
understand the dynamics at characteristic branch points.

\begin{lemma}[Behavior of Arms under 1st Return Map]\label{LEMcyclPermutLocalAndGlobalArms}\lineclear\emph{\cite[Lemma 4.4 \& 4.6]{BS}}
    Let $z$ be a characteristic branch point of period $n$. The first return map $\iter{f}{n}$ of $z$
    either permutes the
    local arms of $z$ transitively or it fixes the local arm pointing towards $c_0$ and permutes the
    remaining ones transitively.

    Let $G$ be a global arm of $z$ which does not contain $c_0$.
    If all local arms are permuted transitively, then either $G$ is mapped homeomorphically by $\iter{f}{n}$
    into another global arm of $z$ so that $c_0\not\in\iter{f}{i}(G)$ for all $i\leq n$ or the associated local
    arm of $G$ is mapped to the local arm pointing towards $c_0$. The local arm pointing to $c_0$ is mapped to the local arm pointing to $c_1$ by $\iter{f}{n}$.  
    
    If the local arm towards $c_0$ is fixed,
    either $G$ is mapped homeomorphically by $\iter{f}{n}$ into another global arm of $z$ with
    $c_0\not\in\iter{f}{i}(G)$ for all $i\leq n$ or the associated local arm of $G$ is mapped to the
    local arm pointing towards $c_1$. \qed
\end{lemma}

We call a characteristic branch point $z$ \emph{evil} if the local arm of $z$
pointing towards $c_0$ is fixed. Otherwise $z$ is called \emph{tame}. Observe that in the tame case
of Lemma~\ref{LEMcyclPermutLocalAndGlobalArms}, a global arm might be mapped homeomorphically without
hitting $c_0$ by the first return map, yet its associated local arm is pointing to $c_0$. While
the term ``evil" is reserved for branch points, we also call inner points whose local arm to $c_0$ is not fixed \emph{tame}. There is a purely combinatorial condition to determine whether a point is tame or not.

\begin{lemma}[Tame Points]\label{LEMcombCondForPermutation}\lineclear
    \emph{\cite[Proposition 4.8]{BS}} Let $z$ be a characteristic
    point of period $n$ and itinerary $\tau$. Then $z$ is tame if and
    only if $n$ is contained in the internal address associated to $\tau$. \qed
\end{lemma}

From this criterion, it follows that the itinerary of an evil branch point is always equal to the
lower kneading sequence $\overline{\upper}(\mu)$ of some $\star$-periodic kneading sequence $\mu$.

\begin{definition}[Admissible Kneading Sequences]\lineclear
    \emph{\cite[Definition 4.9]{BS}} We call a $\star$- or preperiodic kneading sequence $\nu$ \emph{admissible} if its associated Hubbard tree contains no evil  branch point.
    Otherwise $\nu$ is called \emph{non-admissible}.
\end{definition}

\begin{proposition}[Realizable Hubbard Trees]\lineclear\emph{\cite[Proposition 4.10]{BS}}
    A Hubbard tree $(T,f)$ can be embedded into the plane in such a way that $f$ respects the cyclic order of the local arms at all branch points of $T$ if and only if it contains no evil branch point. \qed
\end{proposition}

 From \cite{orsay} and \cite{poirier} it follows that every Hubbard tree which has no evil branch point is realizable by some quadratic polynomial.

For investigating the structure of the set of kneading sequences the following lemma is essential:

\begin{lemma}[Orbit Forcing]\label{LEMorbitForcing} \lineclear\emph{\cite[Lemma 6.2]{BS}}
    Let $(T,f,\nu)$ and $(\widetilde{T},\tilde{f},\tilde{\nu})$ be two Hubbard trees with $\nu$ and $\tilde{\nu}$
    $\star$-periodic. Suppose that one of the following is true:
    \begin{itemize}
        \item There are two characteristic points $p\in T$ and $\tilde{p}\in\widetilde{T}$
            with identical itineraries.
        \item There is a characteristic point $p\in T$ with itinerary $\upper(\tilde{\nu})$
                or $\overline{\upper}(\tilde{\nu})$. We set $\tilde{p}=\tilde{c}_1$.
       \item There is a characteristic point $\tilde{p}\in\widetilde{T}$ with itinerary $\upper(\nu)$ or
                $\overline{\upper}(\nu)$. We set $p=c_1$.
    \end{itemize}
    Then for every characteristic point $p^\prime\in[c_0,p[$ with $\tau(p^\prime)\not=\tau(p)$ \footnote{If the second item holds, we have to exclude the possibility that $\tau(p')=\overline{\upper}(\tilde{\nu})$. For minimal Hubbard trees, which are defined in Section~\ref{SECstructureHT} and which we will work with, this is not necessary.}
 there is a
    characteristic point $\tilde{p}^\prime\in[\tilde{c}_0,\tilde{p}[$ such that $p^\prime$ and
    $\tilde{p}^\prime$ have the same itinerary, the same number of arms and the same type (i.e., tame or
    not). \qed
\end{lemma}
We introduce a partial order $<$ on the set of all $\star$-periodic kneading sequences.
Theorem~\ref{THMkneadingSequAreRealized} and Lemmas \ref{LEMcyclPermutLocalAndGlobalArms} and
\ref{LEMorbitForcing} guarantee that this is well-defined.

\begin{definition}[Order]\label{DEForder}\lineclear
    \emph{\cite[Definition 6.1]{BS}} Let $\nu$ and $\tilde{\nu}$ be two $\star$-periodic kneading sequences. Then
    $ \nu < \tilde{\nu} :\iff $ the Hubbard tree of
                                                $\tilde{\nu}$ contains a characteristic point with itinerary $\upper(\nu)$.
    If $\nu<\tilde{\nu}$ or $\nu=\tilde{\nu}$, we write $\nu\leq\tilde{\nu}$.
\end{definition}

This order extends to all non-periodic sequences $\mu$: we define $\mu>\mu^k$ for all $k$, i.e.\
$\mu$ is larger than all its truncated sequences, and if $\nu$ is $\star$-periodic, then $\nu>\mu$
if and only if $\nu>\mu^k$ for all $k\in\N$. To be able to compare any two sequences in $\Sigma^{\star\star}$, we take
the transitive hull of the order relation defined so far. More details can be found in~\cite[Proposition
6.10]{BS}.

The space $\Sigma^{\star\star}$ together with the order introduced above is called the
\emph{parameter tree}.

\begin{remark}
    In the remainder of this paper, we only investigate Hubbard trees for which the critical point is periodic. So whenever we speak of a Hubbard tree, then we implicitly assume that $c_0$ is periodic.  Note that the associated kneading sequence of such a Hubbard tree
    is $\star$-periodic. This is no restriction when investigating the structure of the space
    of kneading sequences, since any kneading sequence in $\Sigma^{\star\star}$ can be approximated arbitrarily
    well by $\star$-periodic ones.
\end{remark}


\section{Structure of the Hubbard Tree}\label{SECstructureHT}

In this section, we show that for any characteristic point $z$ in a Hubbard tree $(T,f)$
there is a characteristic point $z^\prime\in\, ]z,c_1]$ closest to $z$. If $\tau$ is the
itinerary of $z$ and $n$ its period, then the point $z^\prime$ has itinerary $B_q(\tau)$ or $\upper^{-1}_{qn}(\tau)$ for some $q$. For the
admissible case this means that the arrangement of postcritically finite parameters in
$\mathcal{M}$ is reflected in the arrangement of characteristic points of the Hubbard tree
$T$.

\medskip
We start this section by discussing a way to reduce a Hubbard tree to an equivalent \emph{minimal}
Hubbard tree. For this, we first show some basic properties of Hubbard trees.

\begin{lemma}[Characteristic Point with Minimal Period]\label{LEMexactPeriodPointAndItinerary}\lineclear
    Let $(T,f)$ be a Hubbard tree and $p\in T$ be a characteristic point with itinerary $\tau$.
    Let $m$ be the period of $p$ and $\widetilde{m}$ the period of $\tau$. Then there is a characteristic
    point $\tilde{p}\in T$ that has itinerary $\tau$ and exact period $\widetilde{m}$.
\end{lemma}

\begin{proof} If $m=\widetilde{m}$, then there is nothing to show. Otherwise, observe that $\widetilde{m}$ is a divisor of $m$. Denote by $T_\tau$ the subtree of $T$ that consists of all points with
    itinerary $\tau$. By expansivity, $T_\tau$ is a (not necessarily closed) $n$-od. Moreover, $T_\tau$ contains a point $\tilde{p}$ which has period $\widetilde{m}$, see~\cite[Chapter 3]{BS}. If the branch point $b$ exists, it has period $\widetilde{m}$, and we set $\tilde{p}=b$.

    If $T_\tau$ is a point, we are done since $p=\tilde{p}$.
    In the other two cases we are going to show that $\tilde{p}$ is characteristic. Suppose
    $\tilde{p}$ is not characteristic; then there is an $l<\widetilde{m}$ such that
    $\iter{f}{l}(\tilde{p})\in\, ]\tilde{p},c_1[$.
    If $p\in [c_0,\tilde{p}]$, then $[p,\tilde{p}]\subset \iter{f}{l}([p,\tilde{p}])$. If $\tilde{p}\in[c_0,p]$,
    then we have either $\iter{f}{l}(p)\in\,]\tilde{p},p[$ or $\tilde{p}\in [\iter{f}{l}(p),p]$,
    because $]p,\tilde{p}[$ contains no branch point of $T$. All these possibilities yield
    a periodic point of period (dividing) $l$ by
    the Intermediate Value Theorem. But in $T_\tau$ there is no periodic point of period smaller than $\widetilde{m}$.
\end{proof}

\begin{lemma}[Periods of Periodic Points]\label{LEMperiodsSameItinerary}\lineclear
 Let $(T,f)$ be a Hubbard tree, $\tau\in\{\0,\1\}^{\N}$ be a periodic
itinerary of period $n$ and $T_{\tau}$ the set of all points with itinerary $\tau$. Then there
exists a $k\in \N$ such that the period of any periodic point with itinerary $\tau$ either equals
$n$ or $k\cdot n$.

More precisely, there is at most one periodic point $b\in T_{\tau}$ of period $n$ such that not all
of its local arms are fixed under $\iter{f}{n}$. If such a point $b$ exists, then $k>1$ is the period of any local arm at $b$ which is not fixed. All periodic points in $T_{\tau}$ which are contained in non-fixed arms of
$b$ have period $k\cdot n$ whereas the ones contained in the fixed arm have period $n$. If $b$
does not exist, then all periodic points of itinerary $\tau$ have period $n$.
\end{lemma}

\begin{proof}
Observe that
$\iter{f}{i}|_{T_\tau}$ is a homeomorphism for all $i\in\N$. We know that the existence of a
periodic point with itinerary $\tau$ forces the existence of a periodic point $p$ of period $n$ in
$T_{\tau}$. Let $p'\in T_{\tau}$ be any other periodic point. If $\iter{f}{n}(p')$ is in the same
global arm of $p$ as $p'$, then $\iter{f}{n}(p')=p'$ because otherwise $p'$ would have an infinite
orbit. If $\iter{f}{n}(p')$ is in a different global arm of $p$ than $p'$, consider the first time
$kn$ that the local arm containing $p$ is mapped to itself. Then again $\iter{f}{kn}(p')\not=p'$
would force an infinite orbit for $p'$.

It is easy to see that there is at most one $n$-periodic point in $T_{\tau}$ such that not all of
its local arms are fixed (and that this point is the unique branch point of $T_{\tau}$ if it
exists). Hence the period of any periodic point in $T_{\tau}$ is either $n$ or $kn$ and the rest of
the statement follows from Lemma~\ref{LEMcyclPermutLocalAndGlobalArms}.
\end{proof}

\begin{definition}[Attracting Dynamics]\lineclear
A Hubbard tree $(T,f)$ with $n$-periodic critical point has \emph{attracting dynamics} if for each $c_i\in \orb(c_0)$, there is a neighborhood $U_i$ of $c_i$
 such that for all $x\in U_i$, $\iter{f}{jn}(x)\to c_i$ as $j\to\infty$.
\end{definition}

The following properties of Hubbard trees with attracting dynamics are immediate.
\begin{corollary}[Hubbard Trees with Attracting Dynamics]\label{CORpropAttractingDynamics}\lineclear
Let $(T,f)$ be a Hubbard tree with attracting dynamics. 

\begin{enumerate}
	\item[(i)] If $p$ is a limit point of periodic points $p_n$ which all have the same itinerary $\tau$, then  	$\tau(p)=\tau$ and $p$ is periodic.  

	\item[(ii)]  There is a neighborhood $U$ of $c_1$ such that $\iter{f}{j}|_U$ is a homeomorphism	for all $j\in\N_0$.
\end{enumerate}
\end{corollary}

\begin{proof}
By  Lemma \ref{LEMperiodsSameItinerary}, there is an $m\in\N$ such that all $p_n$ have period $m$. By continuity, $\iter{f}{m}(p)=p$, i.e., $p$ is periodic and its
period divides $m$. Moreover, $p$ has either itinerary $\tau$ or is on the
critical orbit. But since $T$ has locally attracting dynamics at the critical orbit, no point on
the critical orbit can be the limit point of periodic points. 

Let $n$ be the period of $c_1$. There is a neighborhood $U$ of $c_1$ such that $U$ contains no precritical point of step at most $n$. Since $(T,f)$ has attracting dynamics, we can choose $U$ so small that for all $p\in U$, $\iter{f}{n}(p)\in\,]p,c_1[$. If there is a precritical point $\xi\in U$, then there is a $j_0$ such that $\iter{f}{j_0}(\xi)$ is precritical of step at most $n $ and by the choice of $U$, we have $\iter{f}{j_0}(\xi)\in U$, a contradiction.
\end{proof}

\begin{lemma}[Representative with Attracting Dynamics]\label{LEMattractingHT}\lineclear
Every equivalence class of Hubbard trees contains a representative that has attracting dynamics.
\end{lemma}

\begin{proof}
    For any given equivalence class pick a Hubbard tree $(T,f)$ and suppose that the critical orbit is
    not attracting. We are going to define a new dynamics $\tilde{f}$ on $T$ such that any
    point on the critical orbit is locally attracting. To achieve this it suffices to change $f$ locally at $c_0$. Let $V=\{p\in T: p \mbox{ is marked}\}$, $n$ be the period of the critical point $c_0$ and $G$ the global arm of $c_0$ whose associated local arm is fixed under $\iter{f}{n}$.
    Choose $y_0\in G$ such that the interval $I:=\,]c_0,y_0[$ has the following properties: $\iter{f}{n}|_I$ is a homeomorphism onto its image, $I\cap V =\emptyset=\iter{f}{n}(I)\cap V$ and $\iter{f}{i}\cap I=\emptyset$ for all $0<i<n$.
 Without loss of generality, we can assume that there is a $z\in I$ such that $\iter{f}{n}(z)$=z: if such a point does not exist then $p\in\,]c_0,\iter{f}{n}(p)[$ for all $p\in I$ and  we can pick $z,y,y'\in I$ with $c_0<y<z<\iter{f}{n}(z)<y'$ ($<$ denotes the natural order on $\overline{I}$ with $c_0$ as the smallest element). There is a homeomorphism $h:T\to T$ such that $h|_{T\setminus I}\equiv \mbox{id}$ and $h(\iter{f}{n}(z))=z$. Set $f':=h\circ f$. Then $\iter{(f')}{n}(z)=z$ and the two Hubbard trees $(T,f)$ and $(T,f')$ are equivalent.  
 
Now let $z\in I$ be fixed under $\iter{f}{n}$. Pick any homeomorphism $\varphi: [c_0,z]\to[0,1]$ with $\varphi(c_0)=0$ and consider the function $h:[0,1]\to[0,1],\,x\mapsto x^2$. It induces a map 
 \[
 \tilde{h}: T\longrightarrow T,\; p \mapsto \left\{ \begin{array}{ll}\varphi^{-1}\circ h \circ \varphi (p) & \mbox{if  } p\in [c_0,z]\\ p &\mbox{otherwise} \end{array} \right..
 \]
 Let $f_{-n}$ be the inverse branch of $\iter{f}{n}|_I$ that maps $[c_0,z]$ onto itself and define
\[
g: T\longrightarrow T,\; p \mapsto \left\{ \begin{array}{ll}\tilde{h}\circ f_{-n}(p) & \mbox{if  } p\in [c_0,z]\\ p &\mbox{otherwise} \end{array} \right..
\]
 For the continuous map $\tilde{f}:=g\circ f$, $(T,\tilde{f})$ is a Hubbard tree that is equivalent to the given one: for this observe that $V\subset T\setminus f^{-1}(]c_0,z[)=:S$ and that $\tilde{f}|_S=g\circ f|_S=f|_S$. Moreover, $(T,\tilde{f})$ is minimal: it suffices to show that $\iter{f}{jn}(x)\to c_0$ as $j\to\infty$ for all $x\in\, ]c_0,z[$. For any $x\in\,]c_0,z[$, we have 
\[	
\iter{\tilde{f}}{n}(x)=(g\circ f)\circ\iter{(\mbox{id}\circ f)}{n-1}(x)=g\circ \iter{f}{n}(x) = \tilde{h}\circ f_{-n}\circ \iter{f}{n}(x)=\tilde{h}(x)
\]
and $\tilde{h}(x)\in\,]c_0,z[$. Thus, $\iter{f}{jn}(x)=\iter{\tilde{h}}{j}(x)$ for all $j\in\N$ and $x\in\,]c_0,z[$, and $\iter{\tilde{h}}{j}(x)\to c_0$ as $j\to\infty$. 
\end{proof}

If $(T,f)$ is Hubbard tree with attracting dynamics as described in the previous lemma, then
each precritical point has a neighborhood that contains no periodic point.


\begin{definition}[Minimal Hubbard Tree]\lineclear
    We call a Hubbard tree $(T,f)$ \emph{minimal} if it has attracting dynamics and there are no two periodic points sharing the same
    itinerary.
\end{definition}

\begin{remark}
    A minimal Hubbard tree does not contain two distinct preperiodic points of the same itinerary
    either: if there were two such points, say $x$ and $y$, then they would give rise to two distinct periodic points with
    the same itinerary since $\iter{f}{n}|_{[x,y]}$ is a homeomorphism for all $n\in\N$.
    If $(T,f)$ is a minimal Hubbard tree and $p$ is periodic, then the period of $p$ coincides
    with the period of its itinerary by Lemma~\ref{LEMexactPeriodPointAndItinerary}.
\end{remark}

\medskip
The following statements are immediate corollaries to the definition of minimal Hubbard trees. The crucial part is that minimal Hubbard trees have attracting dynamics.

\begin{lemma}[Periodic Points are Repelling]\label{LEMperiodicPtsRepelling}\lineclear
Let $(T,f)$ be a minimal Hubbard tree. If $z\in T$ is an $n$-periodic point disjoint from the critical orbit, then $z$ is repelling, i.e., there is a neighborhood $U$ of $z$ such that for all $p\in U$, $p\in\,]z,\iter{f}{j_pn}(p)[$, where $j_p$ is the period of the local arm at $z$ pointing to $p$.
\end{lemma}

\begin{proof} Let $k$ equal the period of a non-fixed local arm at $z$ if such a local arm exist, and otherwise set $k=1$. Pick $U$ so small that $\iter{f}{kn}|_U$ is a homeomorphism. For any global arm $G$ of $z$ set $L:=G\cap U$. By minimality, either $\iter{f}{jn}(p)\in\,]z,p[$ for all $p\in L$ or $p\in\,]z,\iter{f}{jn}(p)[$ for all $p\in L$, where $j\in\{1,k\}$ is the period of the local arm associated to $G$. Suppose the first case holds.
Then $\iter{f}{jn}(p)\in\,]z,p[$ extends to all $p\in L':=\{x\in G:\,\tau(x)=\tau(z)\}\supset L$. The set $L'$ is an interval and its boundary point $p_0\not=z$ is a precritical point. By continuity, $p_0$ is periodic and thus on the critical orbit, yet all points in $L'$ (except for $z$) are repelled by $p_0$. This contradicts that $(T,f)$ has attracting dynamics. 
\end{proof}

\begin{lemma}[Existence of Dynamical Parent]\label{LEMexistenceDynMother}\lineclear
Let $(T,f,\nu)$ be a minimal Hubbard tree. Then there is a characteristic point $z$ such that $\tau(z)=\overline{\upper}(\nu)$.
\end{lemma}

\begin{proof}
Let $n$ be the period of $\nu$ and let $U\subset T$ be the maximal connected set with $c_1\in\overline{U}$ such that $U$ contains no precritical point and such that for all $p\in U$, $\iter{f}{jn}(p)$ converges to $c_1$ as $j\to\infty$. By Lemma \ref{CORpropAttractingDynamics}, $U\not=\emptyset$ and since each branch point in $T$ has finite orbit, $U$ is an interval contained in $]c_0,c_1[$. Clearly $\iter{f}{n}(U)\subset U$. Since $(T,f)$ has attracting dynamics, $\orb(c_0)\cap \overline{U} =\{c_1\}$. Let $z\in \partial{U}$ and $z\not=c_1$. Then by continuity, $\iter{f}{n}(z)\in[z,c_1]=\overline{U}$, and if $\iter{f}{n}(z)\in\, ]z,c_1[$, then $z$ is not precritical. So again by
continuity, there is a $y\not\in [z,c_1]$ close to $z$ such that $[y,z]$ contains no precritical point of step at most $n$ and $\iter{f}{n}(y)\in\, ]z,c_1[$. Thus $\iter{f}{n}(]y,z[)\subset U$ and $[y,z]$ contains no precritical point. Moreover, $\iter{f}{jn}(y)$ converges to $c_1$. This contradicts that $U$ is maximal. Therefore $\iter{f}{n}(z)=z$. All points in $[z,c_1[$ have the same itinerary, which equals $\overline{\upper}(\nu)$ by \cite[Lemma 5.16]{BS}. This proves the claim.
\end{proof}

\begin{proposition}[Existence of Minimal Trees] \label{PROPminimalHT}\lineclear
    Every equivalence class of Hubbard trees contains a minimal Hubbard tree.
\end{proposition}

\begin{proof} Let $(T,f)$ be a representative that has attracting dynamics.
    For any periodic itinerary $\tau$, let $X_{\tau}\subset T$ be the smallest connected subset of $T$
    which contains all periodic points of itinerary $\tau$. By expansivity, $X_{\tau}$
    contains at most one branch point. By Corollary~\ref{CORpropAttractingDynamics}, the set $X_{\tau}$ is a closed subset of $T$.
    We define the following equivalence relation on $T$:
   	 \[ x\sim y  \quad:\iff
                \left.\begin{array}{ll} \mbox{the itineraries of $x,y$ coincide and }
                \exists \;n\in\N:\;\iter{f}{n}(x) \mbox{ and} \\
          \iter{f}{n}(y)\in X_{\tau} \mbox{ for some periodic itinerary $\tau$.}\end{array}\right.
	\]
    Observe that an equivalence class is either a singleton, of the form $X_{\tau}$ or a connected component of a non-periodic iterated
    preimage of some $X_{\tau}$. Thus, all equivalence classes are closed subsets of $T$. Moreover,
    the equivalence class of any precritical point is trivial and if $x_0\not\in X_{\tau}$ has
    itinerary $\tau$, then there is a neighborhood $U$ of $x_0$ such that the equivalence class of any $x\in U$
    is trivial.
    Let $\tilde{f}$ be the dynamics induced by $f$ on the quotient $\widetilde{T}:=T/_{\sim}$ and $\pi :T \longrightarrow \widetilde{T}$ the
    natural projection map. We are going to show that
    $(\widetilde{T},\tilde{f}\,)$ is a Hubbard tree equivalent to $(T,f)$.

    We first prove that any two points $x\not=y\in \widetilde{T}$ can be separated by a third point:
    there are endpoints $p_x$ and $p_y$ of $\pi^{-1}(x)$, $\pi^{-1}(y)$ respectively such that $]p_x,p_y[$
    does not intersect $\pi^{-1}(x)\cup \pi^{-1}(y)$. If the itineraries of $p_x$ and $p_y$
    are different, then there is a precritical point $\xi \in\, ]p_x,p_y[$.
    If they have equal itineraries, then the equivalence class of at least one of the two points $p_x,p_y$
    is trivial and thus, there is a $\xi\in\,]p_x,p_y[$ whose equivalence class is trivial as well.
    Hence in both cases, the two components $U,U'$ of
    $T\setminus\{\xi\}$ are disjoint open saturated sets, one of which contains $\pi^{-1}(x)$ and the other
    $\pi^{-1}(y)$. Therefore, $\pi(U),\;\pi(U')$ are open in $\widetilde{T}$, disjoint and contain $x$ and
    $y$, respectively. Since $\widetilde{T}\setminus\{\pi(\xi)\}=\pi(U)\cup\pi(U')$, $x$ and $y$ are
    separated by $\pi(\xi)$.
    This implies that $\widetilde{T}$ is metrizable and since all equivalence classes are connected, $\widetilde{T}$
    is a tree (c.f. \cite[9.42, 9.45]{nadler}).

    Observe that $\tilde{f}$ is continuous: $f(X_{\tau})\subset X_{\sigma(\tau)}$ and thus, for any open
    saturated set $U$, $f^{-1}(U)$ is open and saturated. Since $f^{-1}(X_{\sigma(\tau)})$ splits into at
    most two equivalence classes, $\tilde{f}$ is a local homeomorphism and at most $2$-to-$1$.
    Moreover, $(\widetilde{T},\tilde{f})$ meets the expansivity condition. Putting everything together,
    $(\widetilde{T},\tilde{f})$ is a Hubbard tree and since it has the same kneading sequence as $T$,
    it is equivalent to $(T,f)$.
\end{proof}

In the remainder of this paper, we will restrict ourselves to minimal Hubbard trees, i.e., whenever we speak of a Hubbard tree, we mean its minimal representative.

\begin{lemma}[Points with Bifurcation Itinerary]\label{LEMnoPerPoint_BifurcationItinerary}\lineclear
 Let $(T,f,\nu)$ be a Hubbard tree and $z$ be a characteristic point of period $n$. If $z'$ is periodic with itinerary $B_k(\tau(z))$ for some $k\in\N$, then there is an $0\leq j<k$ such that $\iter{f}{jn}(z')$ is the characteristic point of $\orb(z')$. 
  
Moreover, $\iter{f}{jn}(z')\in\,]z,c_1[$ unless $j=0$, $\tau(z)=\upper(\mu)$ and $\tau(z')=\overline{\upper}(\mu)$ for some $\star$-periodic kneading sequence $\mu$.
\end{lemma}

\begin{proof}
Let $\xi$ be the precritical point in $]z,z'[$ of lowest step, which equals $kn$. Assume first that $z'$ is characteristic. If $k>1$ and $z\in\,]z',c_1[$, then the local arm of $z$ pointing to $c_0$ is fixed under $\iter{f}{n}$. On the other hand $\iter{f}{kn}([\xi,z])=[z,c_1]$, a contradiction. If $k=1$ and  $z\in\,]z',c_1[$, it follows that the local arm of $z'$ pointing to $c_1$ is fixed and therefore the local arm of $z'$ pointing to $c_0$ is fixed, too. Thus $\tau(z)$ equals the upper kneading sequence and $\tau(z')$ the lower kneading sequence of a $\star$-periodic kneading sequence $\mu$ as claimed. 
 
Now suppose that $z'$ is not characteristic and let $\iter{f}{i_0}(z')$ be the characteristic point of $\orb(z')$, where $0<i_0<kn$. Since $\tau_i(p)=\tau_i(z)$ for all $0< i\leq kn$ and for all points $p\in[z,z']$, minimality implies that $\orb(z)\cap\,]z',z[\,=\emptyset$ and $\orb(z')\cap\,]z',z[\,=\emptyset$. Thus, $\iter{f}{i_0}(z')\in\,]z,c_1[$. If $i_0\not=jn$ for all $0<j<k$, then $z\in\,]\iter{f}{i_0}(z),\iter{f}{i_0}(z')[$. Therefore $\iter{f}{kn-i_0}(z)\in\iter{f}{kn-i_0}(]\iter{f}{i_0}(z),\iter{f}{i_0}(z')[)=\,]z,c_1]\cup\,]z',c_1]$. Since $z$ is characteristic, it follows that $\iter{f}{kn-i_0}(z)\in\,]z,z'[$, contradicting that $]z,z'[$ contains no iterate of $z$. 
\end{proof}

There are examples where $z\in\,]z',\iter{f}{jn}(z')[$ and where $z'\in\,]z,\iter{f}{jn}(z')[$: for the first possibility consider the Hubbard tree with kneading sequence $\nu=\overline{1010111\star}$ and the $8$-periodic point $z'$ with itinerary $\overline{10101011}$; for the second case, consider the Hubbard tree with kneading sequence $\nu=B_2(\overline{10101110101\star})$ and the periodic point $z'$ with itinerary $\overline{101010101011}$. In both cases, the periodic point $z$ has itinerary $\tau(z)=\overline{10}$.

\begin{theorem}[Existence of Dynamical Bifurcation Point] \label{THMexistenceBifuracation}\lineclear
    Let $(T,f,\nu)$ be a Hubbard tree and $z\in \,]c_0,c_1[$ be a characteristic $n$-periodic point with itinerary  $\tau$.
    Then there is a characteristic point $z^\prime\in\, ]z,c_1]$ closest to $z$. Moreover,
    either $z^\prime\in\,]z,c_1[$ and $z'$ has itinerary $B_Q(\tau)$, or
    $z^\prime=c_1$ and $\nu=\upper^{-1}_{Qn}(\tau)$, where $Q$ is the period of the local arm at $z$ that  points to $c_1$.
\end{theorem}

By Lemma \ref{LEMcyclPermutLocalAndGlobalArms}, $Q$ or $Q+1$ equals the number of arms at $z$according as $z$ is tame or not.

\begin{proof}
Let $G$ be the global arm of $z$ containing the critical value $c_1$ and set $N:=Qn$. Moreover, let $H$ be the connected component of the set $\{ p\in \overline{G}: \iter{f}{N}(p)\in\overline{G} \}$ that contains $z$. This set $H$ contains a non-empty interval $I:=[z,y_0[\,\subset \overline{G}$ such that $\iter{f}{N}|_I$ is a homeomorphism, $\mathring{I}$ contains no marked point and for all $p\in \overline{I}$, $p\in\,]z,\iter{f}{N}(p)[$ (this is possible by Lemma~\ref{LEMperiodicPtsRepelling}).

Let $\xi$ denote the precritical point in $H$ of lowest step. We show first that {\sc step}$(\xi)=N$. Note that there is a $0<j<N$ such that $c_0\in \iter{f}{j}(H)$: if such an iterate does not exist, then $H=\overline{G}$. But then $\iter{f}{N}(G)=G$ and $c_0\not\in\iter{f}{j}(G)$ for all $j\in\N$, which contradicts that $c_0$ is periodic. Therefore, there is a smallest integer $j_0$ with $c_0\in\iter{f}{j_0}(H)$. By  Lemma~\ref{LEMcyclPermutLocalAndGlobalArms}, $j_0\not=kn-1$ for $kn<N$ and thus, if $j_0<N-1$, then $\iter{f}{j_0+1}(H)\ni c_1$ contains $z$ in its interior. Consequently, $\iter{f}{N-j_0-1}(z)\in \iter{f}{N}(H)\subset G$, contradicting that $z$ is characteristic. Hence, $j_0=N-1$.

As a consequence, $\iter{f}{N}|_{H}$ is at most $2$-to-$1$, $\overline{G\setminus H}$ consists of at most one connected component and $H\cap(\overline{G\setminus H})=\{z_{-}\}$ with $\iter{f}{N}(z_{-})=z$. 
Consider the set $H':=H\setminus I$ and the continuous map $g:=r\circ\iter{f}{N}:H'\longrightarrow H'$, where $r:\overline{G}\longrightarrow H'$ is the unique retraction. Since the tree $H'$ has the fixed point property, there is a point $z'\in H'$ with $g(z')=z'$. If $g(z')\not=\iter{f}{N}(z')$, we have that either $z'=y_0$ or $z'=z_{-}$. But by the choice of $I$, $g(y_0)\not=y_0$ and since $z_{-}$ is an $\iter{f}{N}$-preimage of $z$, $g(z_{-})=y_0$. Thus $g(z')=\iter{f}{N}(z')$, which shows the existence of an $N$-periodic point in $H$. By minimality, $\xi\in\,]z,z']$ ($\xi$ as defined above). Consequently, $\tau(z')=B_Q(\tau)$ or $\nu=\upper^{-1}_{Qn}(\tau)$ if $z'=c_1$. 
If $z'\not=c_1$, then Lemmas \ref{LEMnoPerPoint_BifurcationItinerary}  and \ref{LEMcyclPermutLocalAndGlobalArms} imply that $z'\in\,]z,c_1[$ is characteristic.
\end{proof}


\section{Structure of the Parameter Plane}\label{SECparamPlane}
Now we turn to the parameter plane. We start this section by discussing important
consequences of Theorem~\ref{THMexistenceBifuracation}.

\begin{corollary}[Non-Tame Inner and Branch Points]\label{CORevilBranchPointNoUpper}\lineclear
Let $(T,f,\nu)$ be a Hubbard tree and $z\in T$ be a characteristic point of period $n$ and itinerary $\tau$ such that the
first return map fixes its local arm to $c_0$, and set $\mu:=\upper^{-1}_n(\tau)$. Then, if $z$ is
an inner point, $\mu\leq\nu$. If $z$ is a branch point, then $\mu\not<\nu$ but for all
$\star$-periodic $\tilde{\mu}<\mu$ we have that $\tilde{\mu}<\nu$.
\end{corollary}

\begin{proof} 
  If $z$ is an inner point, then by Theorem~\ref{THMexistenceBifuracation}
  either there is a characteristic point $z^\prime\in\,
  ]z,c_1[$ with itinerary $B_1(\tau)=\upper(\mu)$ and $\mu<\nu$, or
  $z^\prime=c_1$ and $\nu=\mu$. 
  If $z$ is a branch point, we have to show that there is no characteristic $n$-periodic point with
  itinerary $\upper(\mu)$.
  By way of contradiction, let us assume that there is such a characteristic point $p\in [c_0,c_1]$.
   By Lemma~\ref{LEMcyclPermutLocalAndGlobalArms}, there are no $n$-periodic points in $[z,c_1]$.
  On the other hand, if $p\in [c_0,z]$, then $c_0\in \iter{f}{(n-1)}([p,z])$ implies that
  $\iter{f}{n}$ maps the local arm at $z$ pointing to $c_0$ to the one pointing towards $c_1$,
  contradicting that $z$ is evil. The last part of the claim is an immediate consequence of
  Orbit Forcing~\ref{LEMorbitForcing}.
\end{proof}

Let us extend the notion of wakes to the set $\Sigma^{\star\star}$: by the \emph{wake} of a kneading sequence $\mu$, we mean the set $\{\mu'\in\Sigma^{\star\star}: \mu\leq\mu'\}$. For any $q\geq 2$  the set  $\{\mu'\in\Sigma^{\star\star}: B_q(\mu)\leq\mu'\}$ is called a \emph{subwake} of $\mu$. 

Now suppose that $z$ is an evil branch point of the Hubbard tree $(T,f,\nu)$. Let $q+1$ be the number of its arms, $n$ its period, and let $\mu$ be $\star$-periodic such that $\tau(z)=\overline{\upper}(\mu)$. Then either there
is a next characteristic point $z^\prime$ after $z$, which has itinerary
$B_{q}(\overline{\upper}(\mu))$, or $z^\prime=c_1$ and $\nu=\overline{B}_{q}(\mu)$. So in a way,
$\nu$ is contained in the wake of a bifurcation sequence of $\mu$, only that we used the lower
kneading sequence instead of the upper one to build the bifurcation sequence. Contrary to
``normal" bifurcations, we do not have that $\mu<\nu$. To distinguish this kind of bifurcation from
the classical one, we call it a \emph{backward bifurcation}. The set  $\{\mu'\in\Sigma^{\star\star}: \overline{B}_{q}(\mu)\leq \mu'\}$ is called a \emph{non-admissible subwake} of $\mu$. Backward bifurcations can only happen
at a $\star$-periodic kneading sequence $\mu$ whose lower kneading sequence has the same period as
$\mu$ itself: a backward bifurcation needs the existence of an evil branch point with itinerary
$\overline{\upper}(\mu)$, and the period of a branch point coincides with the period of its
itinerary. For
admissible kneading sequences this means that backward bifurcations, and thus branching off into non-admissibility, only occurs at primitive
hyperbolic components: for all bifurcation components, the associated lower kneading sequence has
period strictly dividing the period of the component itself (c.f.~\cite{lausch}).

\begin{corollary}[No Irrational Decorations]\label{CORq-wake}\lineclear
   Let $\nu\not=\nu'$ be two $\star$-periodic kneading sequences.
   If $\nu<\nu'$,
   then there is a $q\in\N$ such that $B_q(\nu)\leq\nu'$, i.e., $\nu'$ lies
  in the wake of a bifurcation sequence of $\nu$. In particular, if $\nu$ and $\nu'$ are admissible, then $\nu'$ lies in a $\frac{p}{q}$-subwake of $\nu$.

   If $\nu\not<\nu'$ but $\mu<\nu$ for all $\star$-periodic $\mu<\nu$, then there is a $q$ such that
  $ \overline{B}_q(\nu)\leq\nu'$, i.e., $\nu'$ lies in the wake of a backward bifurcation sequence of $\nu$.
\qed
\end{corollary}

By the definition of the order on $\Sigma^{\star\star}$, the above statement is also true for preperiodic kneading sequences. Using Theorem \ref{THMexistenceBifuracation} and the notion of wakes, we get the following description of the location of all non-admissible kneading sequences. This strengthens the answer given in \cite{BS} to Kauko's question what kind of non-admissible kneading sequences there are  \cite{kauko2}.

\begin{corollary}[Location of Non-Admissible Sequences]\lineclear
 Every non-admissible kneading sequence $\nu$ is contained in the wake of a backward bifurcation sequence of a
 primitive $\star$-periodic kneading sequence that is admissible. \qed
\end{corollary}

\subsection{Orbit Forcing}

Many of the subsequent proofs will be based on iterating triods in $T$ ho\-meo\-mor\-phi\-cally. A \emph{triod} is the convex hull $[x,y,z]$ of three pairwise distinct points $x,y,z \in T$. We call these three points the \emph{generating} points of the triod. A
\emph{non-degenerate triod} is a metric space homeomorphic to the letter $Y$. A \emph{degenerate
triod} is homeomorphic to a closed interval. Observe that a triod $[x,y,z]$ can be pushed forward homeomorphically if and only if $c_0$ is not contained in
its interior. If it is and we want to push forward $[x,y,z]$ homeomorphically, we have to \emph{chop} this triod before. The
chopping must happen in such a way that the resulting triod does not longer contain $c_0$ in its interior and is topologically the same as the
original one.  More precisely, we require that it contains two points out of $\{f(x),f(y),f(z)\}$ and the third one is \emph{chopped off}, that is, replaced by a point $p$ distinct from the two not chopped points. The mutual location of the new generating points must be the same as the one of $x,y,z$. Usually, it suffices to choose $p=c_0$. However, sometimes we want the endpoints to have specific itineraries and thus we have to replace the separated point in a more tricky way.
Since $c_0$ is not a branch point, the procedure of chopping a triod is not possible if and only if
one of the generating points is mapped onto $c_0$ and the other two are contained in two different
components of $T\setminus\{c_0\}$. This event is called \emph{stop case} because it prevents any
further iteration. The \emph{stop case} can only occur for degenerate triods. In the non-degenerate case, we always chop off the generating point which is separated from the other two by the critical point $c_0$. 

\medskip
The content of the following lemma is close to Orbit Forcing~\ref{LEMorbitForcing}: in some
sense it is stronger since we drop the requirement that the point corresponding to the kneading
sequence in whose Hubbard tree we want to force another point must be characteristic. It is an
analogue of the correspondence of dynamic and parameter rays for quadratic polynomials
(c.f.~\cite{lavaurs, milnor}).

\begin{lemma}[Forcing of Characteristic Points]\label{LEMcharPointsAreForced}\lineclear
    Let $(T,f)$ be a Hubbard tree and $p\in T$ be a periodic non-precritical point of period
    $n$ and itinerary $\tau$ such that  $]c_0,c_1] \cap \,]c_0,p]=:\,]c_0,b]$ is not empty. Set $\tilde{\nu}:=\mathcal{A}_n^{-1}(\tau)$ and let $(\widetilde{T},\tilde{f})$ be the Hubbard tree
    associated to the kneading sequence $\tilde{\nu}$.
    If $z\in\, ]c_0,b[$ is a characteristic point,
    then there is a characteristic point $\tilde{z}\in\widetilde{T}$ 
    such that $z$ and $\tilde{z}$ have the same itinerary and are of the same type. Moreover, if $]z,b[$ contains a characteristic point, then $z$ and $\tilde{z}$ have  the same number of arms.
\end{lemma}

\begin{proof}
    The Hubbard tree $(\widetilde{T},\tilde{f})$ associated to $\tilde{\nu}$ exists by Theorem \ref{THMkneadingSequAreRealized}. In order to tell points in the two Hubbard trees
    apart, all points in $\widetilde{T}$ are marked by $\sim$. 
    We can assume that $T$ contains both preimages $p_0,-p_0$ of $p$: if it does not contain the preimage, say, $p_0$,
    then we attach an arc $[c_0,p_0]$ to the tree $T$ such that $[c_0,p_0]$ is mapped homeomorphically onto $[c_1,p]$. This extended tree is not a Hubbard tree in the strict sense anymore because not all of its endpoints are on the critical orbit. However, all other properties of Hubbard trees are preserved under this extension.
    
    Let $p_0$ be the preimage which is not separated from $p$ by $c_0$ and $b_0$ such that $[c_0,b_0]=[c_0,c_1]\cap[c_0,p_0]$. Let $G$ be the global arm of $b_0$ that contains the critical value $c_1$. If $p\not\in G$ then there is no characteristic point in $]c_0,b[$ and the statement is empty.  So from now on, assume that $p\in G$.
    Figure \ref{FigObitForcing} illustrates the location of
    the mentioned points in the tree $T$.

    \begin{figure}[htb] \setlength{\unitlength}{1cm}
    \begin{center}  \begin{picture}(7,4)
                        \epsfig{file=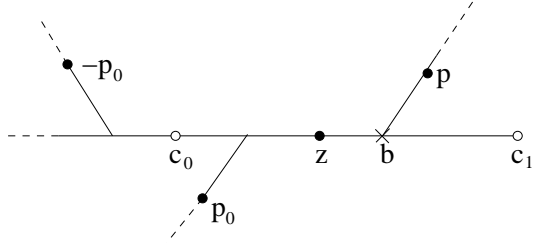, width=7cm}
                    \end{picture}
     \end{center}
     \caption{The point $p$ and its preimages in the Hubbard tree $T$.} \label{FigObitForcing}
    \end{figure}

   Observe that $\orb(z)\cap[p_0,-p_0]=\emptyset$ because $z$ is characteristic. We are going to construct closed intervals $I_k\subset T$ with endpoints in
\[    	P=\{-p_0,p_0,p,f(p),\ldots,\iter{f}{(n-1)}(p)\} \]
    such that $\iter{f}{k}(z)\in I_k$ and $f|_{I_k}$ is a homeomorphism.
     We define the intervals $I_k$ iteratively. Let $m$ be the period of $z$.\lineclear
    $k=0$: Since $z\not\in [p_0,-p_0]$, the interval $I_0:=[p_0,p]$ contains $z$ and clearly $f|_{[p_0,p]}$ is a homeomorphism.\lineclear
    $k\rightarrow k+1$: Suppose $I_k=[x,y]$ has already been defined. By definition, $f$ maps $[x,y]$ homeomorphically onto $[f(x),f(y)]$. If $c_0\not\in[f(x),f(y)]$, then $f|_{[f(x),f(y)]}$ is a homeomorphism, and we set $I_{k+1}:=[f(x),f(y)]\ni \iter{f}{k+1}(z)$. 
    If $c_0\in[f(x),f(y)]$, then $\iter{f}{k+1}(z)$ is either contained in $[\pm p_0,f(x)]$ or in $[\pm p_0,f(y)]$, where $\pm p_0$ denotes the preimage of $p$ such that $c_0$ is not contained in the respective  interval. Define $I_{k+1}$ to be the one of these two intervals that contains $\iter{f}{k+1}(z)$. By definition, $f|_{I_{k+1}}$ is a homeomorphism and since $\orb(z)\cap[p_0,-p_0]=\emptyset$, $\iter{f}{k+1}(z)\in I_{k+1}$. (Note that in general $I_{k+1}\not\subset f(I_k)$.)
    
    This way, we define intervals until $I_{jm}=I_{j'm}$ for some $0\leq j<j'$. Such $j,j'$ exist because $P$ has only finitely many elements. For all $k>j'm$, we set $I_k=I_{k\mod(j'-j)m}$. 
   
    Now we are going to define analogous closed intervals $\widetilde{I}_k$ in the Hubbard tree $\widetilde{T}$
    associated to the kneading sequence
    $\tilde{\nu}$.
    The critical value $\tilde{c}_1=:\tilde{p}$ has itinerary $\tilde{\nu}$.
    We have that  $-\tilde{p}_0=\tilde{p}_0=\tilde{c}_0$. For all $k\in\N_0$, let $\widetilde{I}_k$ be the
    closed arc with endpoints in $ \{-\tilde{p}_0,\tilde{p}_0,\tilde{p},\ldots,\iter{f}{(n-1)}(\tilde{p})\}$
    corresponding to the ones of $I_{k+jm}$.
    Since $p$ and $\tilde{p}$ have the same itinerary (except for the $\star$), $\widetilde{I}_k$ maps homeomorphically to an
    arc containing $\widetilde{I}_{k+1}$. Consider the set
      	\[ 
      	\widetilde{S}_l:=\{\tilde{x}\in [\tilde{c}_0,\tilde{c}_1]: \, \tilde{f}^{\circ k}(\tilde{x}) \in \widetilde{I}_k \text{ for all } k<l\}.
	\]
    $\widetilde{S}_{l+1}\subset \widetilde{S}_l$ and for each $l\in\N_0$, the set $\widetilde{S}_l$ is a non-empty, compact
    interval. Hence $\widetilde{S}:=\bigcap \widetilde{S}_l$
    is also non-empty, compact and connected and each $x$ in the interior of $\widetilde{S}$ has itinerary $\tau (z)$. By definition, $\iter{\tilde{f}}{(j'-j)m}(\widetilde{S})=\widetilde{S}$. Hence there is a periodic point
    $\tilde{z}\in\widetilde{S}$ of period dividing $(j'-j)m$. Suppose there is no periodic point in $\widetilde{S}$ with itinerary $\tau(z)$. Then there are no periodic points in the interior of  $\widetilde{S}$.
     If $\widetilde{S}$ is a non-degenerate interval, then the two endpoints must be (pre-)critical and fixed by $\iter{f}{(j'-j)m}$, that is, they must actually lie on the critical orbit. Thus both endpoints are attracting and there must be a periodic point in the interior of $\widetilde{S}$, a contradiction. This implies that $\widetilde{S}=\{\tilde{x}\}$ and $\tilde{x}\in\orb(\tilde{c}_1)$. The way the set $\widetilde{S}$ was constructed, it follows that $\tilde{x}$ is the limit of precritical points. But this contradicts that $(\widetilde{T},\tilde{f})$ has attracting dynamics.
    Thus, we showed the existence of a periodic point $\tilde{z}$ of itinerary $\tau(z)$ and by minimality, the period of $\tilde{z}$ is $m$.
    
    It remains to show that $\tilde{z}$ is characteristic: \label{PAGEzTildeChar} suppose that it was not, then there is an $l\leq m$ such that
    $\iter{f}{(l-1)}(\tilde{z})=:\tilde{z}_l\in\, ]\tilde{z},\tilde{p}[$. Hence the three points
    $\tilde{z}_l,\,\tilde{z}$ and $\tilde{p}$ form a degenerate
    triod $\widetilde{Y}$ with $\tilde{z}_l$ in the middle ($\tilde{p}=\tilde{c}_1$). But in the Hubbard tree $T$, the points
    $z_l,\,z,\,p$ form a degenerate triod $Y$ with $z$ in the middle because $z$ was characteristic.
    Since $m$ is the exact period of
    $\tau(z)$, we have that $\tau(z)\not=\tau(z_l)$. Therefore, there
    is a first time $k<m$ such that $c_0\in \iter{f}{k}([z,z_l])$.
    Iterate $Y$ until $c_0\in \iter{f}{j}(Y)$ for the first time. If $j< k$, then $c_0\in \iter{f}{j}([z,p])$.
    In this case, recall that no image of $z$ is contained in $[-p_0,p_0]$. Hence, we can replace
    $\iter{f}{j}(p)$ by the one point of $\{-p_0,\,p_0\}$ that is not separated from $\iter{f}{j}(z)$ by $c_0$ without changing the mutual location of the generating points.
    In the Hubbard tree $\widetilde{T}$, we replace $\iter{f}{j}(\tilde{p})$ by $\tilde{c}_0$.
    This ensures that corresponding endpoints of the triods in $T$ and $\widetilde{T}$ still have the same itineraries.
    We go on pushing the modified triod forward until we reach time $k$. This might include
    further choppings of the just described kind.
    Now $c_0\in \iter{f}{k}([z_l,z])$ and the arrangement of points in $T$ implies that the $k$-th image of
    $z$ and $p$ (or of the point which replaced $p$) are on the same side of $c_0$, whereas in $\widetilde{T}$ the
    $k$-th images of $\tilde{z}$ and $\tilde{p}$ (or of the point which replaced $\tilde{p}$) are
    separated by $\tilde{c}_0$, contradicting that $z$ and $\tilde{z}$ as well as $p$ and $\tilde{p}$ (or the points that replaced them)
    have the same itinerary (up to the symbol $\star$).

    The two points $z$ and $\tilde{z}$ are of
    the same type because  the type is completely encoded into the internal address by
    Lemma~\ref{LEMcombCondForPermutation}.
            
    Now let $z$ be a characteristic $m$-periodic point such that $]z,b[$ contains a further characteristic point and let $q$,
    $\tilde{q}$ be the number of arms at $z$, $\tilde{z}$, respectively. We are going to show that $q=\tilde{q}$.
    By Theorem~\ref{THMexistenceBifuracation}, there are characteristic points $z^q\in\,]z,b[\subset T$
    and $\tilde{z}^{\tilde{q}}\in \widetilde{T}$ with itinerary $B_q(\tau(z))$ and
    $B_{\tilde{q}}(\tau(z))$ (or $\upper^{-1}_{\tilde{q}m}(\tau(z))$), respectively. If $q \not=\tilde{q}$, then there is a
    characteristic point $\tilde{z}^q\in\widetilde{T}$ with itinerary $B_q(\tau(z))$  as we just have proven and $\tilde{z}^q\in\,]\tilde{z},\tilde{c}_1[$. Let us assume that
    $\tilde{z}^q\in\,]\tilde{z},\tilde{z}^{\tilde{q}}[$, the case that $\tilde{z}^{\tilde{q}}\in\,]\tilde{z},\tilde{z}^q[$ works exactly the same way. The precritical point $\xi\in[\tilde{z},\tilde{z}^q]$ of smallest step
    has {\sc step}$(\xi)=qm$ and the one in $[\tilde{z},\tilde{z}^{\tilde{q}}]$, denoted by
    $\tilde{\xi}$, has {\sc step}$(\tilde{\xi})=\tilde{q}m$. Therefore we must have that $q>\tilde{q}$ and $\tilde{\xi}\in\,]\tilde{z}^q,\tilde{z}^{\tilde{q}}[$. By minimality, we have for all $\tilde{x}\in\,]\tilde{z},\tilde{\xi}[$ that $\tilde{x}\in\,]\tilde{z},\iter{f}{\tilde{q}m}(\tilde{x})[$, in particular this holds for $\tilde{z}^q$, which consequently cannot be characteristic.
\end{proof}

\begin{remark} 
Let $\mu$ be $\star$-periodic. If $T$ contains no periodic point $p$ with itinerary $\upper(\mu)$, then we can extend the Hubbard tree $(T,f)$ such as to contain
the orbit of $p$. One can check that this extended tree $(T^\prime,f^\prime)$
satisfies all conditions of a Hubbard tree except that not all endpoints are on the critical orbit.
There are several possibilities for the location of $p$ in $T^\prime$:

If $p$ is an element of the component of $T^\prime\setminus\{c_0\}$ not containing $c_1$, then the statement is empty. If $p$ is contained in a subtree branching off from $[c_0,c_1]$,
then we have exactly the same situation as in Lemma~\ref{LEMcharPointsAreForced}.
The last possibility is that $c_1\in\,]c_0,p[$. In this case, $]-p_0,p_0[\,\cap \,T=\{c_0\}$ and $c_0$ is a branch point of $T'$. For all characteristic points $z\in\,]c_0,c_1[$ we can define intervals $I_k$ such that for all $k\in\N_0$, $\iter{f}{k}(z)\in I_k$, $f|_{I_k}$ is a homeomorphism and $\partial I_k\in \orb(p)\cup\{-p_0,p_0\}$. We modify the interval $f(I_k)$ in order to get $I_{k+1}$ if and only if $f|_{I_k}$ is not a homeomorphism. In general, the intervals $I_k$ may contain the critical point $c_0$ (e.g.\ $I_0$ always does). We construct intervals $\tilde{I}_k$ in $\widetilde{T}$ just the same way as we did in the proof of Lemma \ref{LEMcharPointsAreForced}. In $\widetilde{T}$, we have that $\tilde{f}(\tilde{I}_k)\supset \tilde{I}_{k+1}$ for all $k\in\N_0$.  Thus, the remainder of the above proof carries over, and the statement of  Lemma \ref{LEMcharPointsAreForced} extends to the case that $c_1\in\,]c_0,p[$.
\end{remark}

\begin{lemma}[Not Forced Characteristic Points]\label{LEMcharPointAfterBranchPointNotForced}\lineclear
 Let $(T,f,\nu)$ be a Hubbard Tree
and define $p,b\in T$ and $\tilde{\nu}$ as in Lemma~\ref{LEMcharPointsAreForced}. Suppose that
$z\in\,]b,c_1[$ is a characteristic point such that
either there is a characteristic point $y\not=p$ in $[b,z[$ or $b$ is the limit of characteristic
points in $]c_0,b[$. Then there is no characteristic point $\tilde{z}$ with itinerary $\tau(z)$ in the
Hubbard tree $\widetilde{T}$ of $\tilde{\nu}$.
\end{lemma}

\begin{proof} We state the proof for the case that there is a
characteristic point $y\in[b,z[$. For the case that such a point does not exist but $b$ is the limit
of characteristic points we refer to the footnotes. By way of contradiction we assume
that there is a characteristic point $\tilde{z}\in \widetilde{T}$ such that $\tau(z)=\tau(\tilde{z})$.

The triod $[p,y,z]\subset T$, which is well-defined by hypothesis, is degenerate and has $y$ as an inner point. The point
$\tilde{z}\in \widetilde{T}$ forces a characteristic point $\tilde{y}\in \widetilde{T}$ with
itinerary $\tau(\tilde{y})=\tau(y)$ by Orbit Forcing~\ref{LEMorbitForcing}. Thus there is a
degenerate triod $[\tilde{y},\tilde{z},\tilde{p}]$ with $\tilde{z}$ in the middle (where $\tilde{p}=\tilde{c}_1$).\footnote{We set
$y=b$. Any
    characteristic point in $[c_0,b[\,\subset T$ is forced in $\widetilde{T}$ by Lemma~\ref{LEMcharPointsAreForced}.
    Hence there is a limit point $\tilde{b}\in \widetilde{T}$, which is either the critical value or
    has the same itinerary as $b$ (c.f. Lemma~\ref{LEMlimitOfItineraries}). Because of the existence of $\tilde{z}$,
    the second case must hold. We consider the degenerate triod $[\tilde{b},\tilde{z},\tilde{p}]$. }
First note that the \emph{stop case} never occurs in any of the two triods, so that they can be
iterated forever. There is a smallest number $k$ such that $c_0\in\iter{f}{k}([y,z])$. If it happens that
$c_0\in\iter{f}{l}([p,y])$ for some $l<k$, then replace $\iter{f}{l}(p)$ by $\pm z_0$ (by $\pm z_0$
we mean the appropriate preimage of $z$), which is possible because no image of the
characteristic point $y$ is contained in $[-z_0,z_0]$. This yields a degenerate triod $[\pm
z_0,\iter{f}{l}(y),\iter{f}{l}(z)]$ with $\iter{f}{l}(y)$ in the middle.\footnote{Since $b$
is the limit of characteristic points it follows
    that no image of $b$ is contained in $[-z_0,z_0]$.}
Observe that this also implies $\iter{f}{l}(z)\not=\pm z_0$ and hence
$\iter{f}{l}(\tilde{z})\not=\pm \tilde{z}_0$. Therefore, replacing $\iter{f}{l}(\tilde{p})$ by
$\pm\tilde{z}_0$ yields the degenerate triod
$[\iter{f}{l}(\tilde{y}),\iter{f}{l}(\tilde{z}),\pm \tilde{z}_0]$. We iterate the modified
triods until we reach the time $k$ (which might include further choppings of the just described kind). Then $\iter{f}{k}(y)$ and $\iter{f}{k}(p)$ (or the respective image of the point which replaces $p$) are
on the same side of $c_0$ whereas in $\widetilde{T}$ the images of $\tilde{z}$ and $\tilde{p}$ (or the
point which replaces $\tilde{p}$) are on the same side of $\tilde{c}_0$, a contradiction. 
\end{proof}

\subsection{The Branch Theorem}

For the proof of the Branch Theorem we need the following two technical lemmas.

\begin{lemma}[Points of Equal Itinerary]\label{LEMptsOfSameItinerary}\lineclear
    Let $b$ be a periodic point of the Hubbard tree $(T,f)$. Fix any global arm $G$ of $b$ and let
    $I:=\{p\in G:\, p \text{ has the same
    itinerary as } b\}$. Then $I=[b,x[$, where $x$ is a precritical point.
\end{lemma}

\begin{proof}
By expansivity, $I$ is an interval. If $I=[b,x[$, then $x$ is precritical
because the itineraries of $x$ and $b$ are different and there is no precritical point in $[b,x[$.

If $I=[b,x]$, we show that then $x$ must be periodic, which is impossible by minimality. Since all points in
$[b,x]$ have the same itinerary, we have in particular that $c_0\not\in \iter{f}{i}([b,x])$ for all
$i\in\N$. By Lemma~\ref{LEMcyclPermutLocalAndGlobalArms}, there is an $N$ such that
$\iter{f}{N}(b)=b$ and $\iter{f}{N}([b,x])\subset G$. Since $[b,x]$ is the maximal set of points in $G$ with itinerary $\tau(b)$, continuity of $f$ implies that $\iter{f}{N}([b,x])=[b,x]$ and thus $\iter{f}{N}(x)=x$. 
\end{proof}

\begin{lemma}[Itineraries of Limit Points]\label{LEMlimitOfItineraries}\lineclear
Let $(T,f)$ be a Hubbard tree, $x\in T$ be a point which is not precritical and
$\{x_k\}_{k\in\N}$ be a sequence of points converging to $x$. If $\tau$ is the
itinerary of $x$ and $\tau(x_k)$ the one of $x_k$, then $\tau(x_k)\to \tau$ as $k\to\infty$.
\end{lemma}

\begin{proof}Without loss of generality, we can assume that $x_k<x_{k+1}$ for all $k$, where $<$
    refers to the natural order on the arc $[c_0,x]$ with $c_0$ as the smallest element.
    Since there are only finitely many precritical points of a fixed step, the sequence $\tau(x_k)$
    converges to say $\mu$. 
    By possibly taking a subsequence, we can assume that for
    all $k\geq k_0$ the first $k_0$ entries of $\tau(x_k)$ and $\mu$
    coincide. Fix any $m\in\N$ and let $\xi\in[c_0,x]$ be the closest precritical point to $x$ which has {\sc step}$(\xi)\leq
    m+1$. Since $x$ is not precritical the itineraries of all points in $]\xi,x]$ have the same first $m$ entries.
    There is an $M>m$ such that $x_j\in\, ]\xi,x]$ for all $j>M$ and therefore, we have
    $\tau_m(x_j)=\tau_m$ for all $j>M$ (where $\tau_m$ denotes the $m$-th entry of the sequence $\tau$).
    On the other hand, $\tau_m(x_j)=\mu_m$ for all $j>M$ and thus $\mu_m=\tau_m$.
\end{proof}

    For any Hubbard tree $(T,f,\nu)$, $p_k \to c_1$
    implies that $\tau(p_k)\to \overline{\upper}(\nu)$ as $k\to\infty$ \cite[Lemma 5.16]{BS}.
 
\begin{lemma}[Periodic Points Behind $c_1$]\label{LEMPeriodicPointsBehindc1}\lineclear
Suppose that $(T,f,\nu)$ and $(\widetilde{T},\tilde{f},\tilde{\nu})$ are two Hubbard trees such that $\nu\not=\tilde{\nu}$. If necessary extend $T$ such as to contain a periodic point $p_{\tilde{\nu}}$ with itinerary $\upper(\tilde{\nu})$ and $\widetilde{T}$ such as to contain a periodic point $\tilde{p}_\nu$ with $\tau(\tilde{p}_\nu)=\upper(\nu)$. Then it is not possible that $c_1\in\,]c_0,p_{\tilde{\nu}}[$ and $\tilde{c}_1\in\,]\tilde{c}_0,\tilde{p}_\nu[$.
\end{lemma}

\begin{proof}
By way of contradiction assume that for the extended trees $T'$ and $\widetilde{T}'$ we have that $c_1\in\,]c_0,p_{\tilde{\nu}}[$ and $\tilde{c}_1\in\,]\tilde{c}_0,\tilde{p}_\nu[$. By Lemma \ref{LEMexistenceDynMother}, there is a characteristic point $z\in\,]c_0,c_1[\,\subset T\subset T'$ with itinerary $\overline{\upper}(\nu)$ and by the remark after Lemma \ref{LEMcharPointsAreForced}, there is a characteristic point $\tilde{z}\in\,]\tilde{c}_0,\tilde{c}_1[\,\subset \widetilde{T}'$ with $\tau(\tilde{z})=\overline{\upper}(\nu)$.
Thus the precritical point $\tilde{\xi}$ of lowest step in $]\tilde{z},\tilde{p}_\nu[\,\subset \widetilde{T}'$ has step $n$, if $n$ equals the period of $\nu$. If $\tilde{\xi}\in\,]\tilde{z},\tilde{c_1}[$, then $\iter{f}{n}(c_1)\in\,]\tilde{c}_1,\tilde{p}_\nu[$, which contradicts that $\tilde{c}_1$ is an endpoint of $\widetilde{T}$. If $\tilde{\xi}\in\,]\tilde{c}_1,\tilde{p}_\nu[$, then $\orb(\tilde{c}_1)=\infty$, which is not possible either. Consequently $\tilde{\xi}=\tilde{c}_1$, and $\tilde{\nu}=\tau(\tilde{c}_1)=\nu$, a contradiction.
\end{proof}

\begin{definition}[Combinatorial Arc]\lineclear
For two pre- or $\star$-periodic kneading sequences $\nu\leq\nu^\prime$, we define
$[\nu,\nu^\prime]:=\{\mu\in\Sigma^\star:\, \mu \mbox{ pre- or $\star$-periodic such that }
\nu\leq\mu\leq\nu^\prime \}.$ We call $[\nu,\nu^\prime]$ the \emph{combinatorial arc} between $\nu$
and $\nu^\prime$.
\end{definition}

\begin{thm}{Theorem}{1.1}{Branch Theorem}\lineclear
    Let $\nu\not=\tilde{\nu}$ be $\star$-periodic or preperiodic kneading sequences.
    Then there is a unique kneading sequence $\mu$ such that exactly one of the following cases
    holds:
    \begin{enumerate}
        \item   $[\overline{\star},\nu]\cap[\overline{\star},\tilde{\nu}]=[\overline{\star},\mu]$,
                        where $\mu$ is either $\star$-periodic or preperiodic.
        \item    $[\overline{\star},\nu]\cap[\overline{\star},\tilde{\nu}]=[\overline{\star},\mu]\setminus\{\mu\}$,
                    where $\mu$ is a primitive $\star$-periodic kneading sequence.
    \end{enumerate}
\end{thm}

We can split the above two possibilities into subcases so that for any two $\star$-periodic or
preperiodic kneading sequences $\nu,\,\tilde{\nu}$ exactly one of the following cases holds:
    \begin{enumerate}\renewcommand{\labelenumi}{(\roman{enumi})}
        \item $\nu$ and $\nu'$ can be compared, i.e., we have that either $\nu=\tilde{\nu}$, $\nu<\tilde{\nu}$ or $\tilde{\nu}<\nu$.

        \item There is a $\star$-periodic or preperiodic kneading sequence $\mu$ such that
                $\mu<\nu,\,\mu<\tilde{\nu}$ and if $\tilde{\mu}$ is another $\star$-periodic
                kneading sequence with the same property, then $\tilde{\mu}<\mu$.
        \item \begin{enumerate}
                    \item There is a $\star$-periodic kneading sequence $\mu$ such that $\mu\leq\tilde{\nu}$ and
                            $\mu\not<\nu$ but for all $\tilde{\mu}<\mu$:  $\tilde{\mu}<\nu$.
                    \item There is a $\star$-periodic kneading sequence $\mu$ such that $\mu\leq\nu$ and
                            $\mu\not<\tilde{\nu}$ but for all $\tilde{\mu}<\mu$:  $\tilde{\mu}<\tilde{\nu}$.
                    \item There is a $\star$-periodic kneading sequence $\mu$ such
                            that $\mu\not<\nu,\tilde{\nu}$ but for all $\tilde{\mu}<\mu$:
                            $\tilde{\mu}<\nu,\tilde{\nu}$ and there is no $\tilde{\mu}\not<\mu$ such that
                            $\tilde{\mu}<\nu,\tilde{\nu}$.
                \end{enumerate}
\end{enumerate}
Observe that the cases in $(iii)$ can only occur if at least one of the two given kneading sequences
$\nu$ and $\tilde{\nu}$ is non-admissible.

\begin{proof}
We have that $[\overline{\star},\nu]\cap[\overline{\star},\tilde{\nu}]=\{\mu' \mbox{ $\star$- or preperiodic: } \mu'\leq\nu \mbox{ and } \mu'\leq\tilde{\nu}\}$. This together with the way the order on $\Sigma^{\star\star}$ was defined yields that finding the kneading sequence $\mu$ of the Branch Theorem is equivalent to finding the supremum of the set $\{\mu' \mbox{ $\star$-periodic: } \mu'\leq\nu \mbox{ and } \mu'\leq\tilde{\nu}\}$.
     By the definition of $<$, we can assume that both $\nu$ and $\tilde{\nu}$ are $\star$-periodic.
    Let $T$ be the Hubbard tree of $\nu$ and $\widetilde{T}$ the one of $\tilde{\nu}$. If $T$ (or $\widetilde{T}$) contains a characteristic point with itinerary $\upper(\tilde{\nu})$ (or $\upper(\nu)$), then $\nu>\tilde{\nu}$ (or $\tilde{\nu}>\nu$). This is in particular true if one of the given kneading sequences equals $\overline{\star}$.
    For the remaining cases note that, by possibly enlarging $T$, the tree $T$ contains a periodic point $p$ with itinerary $\upper(\tilde{\nu})$ and by  Lemma \ref{LEMPeriodicPointsBehindc1}, $c_1\not\in[c_0,p]$  (interchange $T$ and $\widetilde{T}$ if necessary). 
    Let $b\in T$ such that $[c_0,p] \cap [c_0,c_1]=[c_0,b]$. The point $b$ is either a branch point or equals $p$.
    Let
       	\[
       		P:=\{z\in [c_0,b[\,: z \text{ is a characteristic point} \}.
	\]
    Suppose that there is a kneading sequence $\nu'=\overline{\1\ldots\1\star}$ such that $\nu^\prime<\nu$ and $\nu^\prime<\tilde{\nu}$ (otherwise $[\overline{\star},\nu]\cap[\overline{\star},\tilde{\nu}]=\{\overline{\star}\}$).
    Then by Lemma \ref{LEMcharPointAfterBranchPointNotForced}, $P$ contains the $\alpha$-fixed point and hence is not empty. Set $a:=\sup(P)$. We distinguish
    two cases:

\begin{enumerate}
     \item The first case is that $a=b$. Observe that in this case $b\not=p$: the local arm of $b$ pointing towards $c_0$ is fixed and all iterates of $b$ are contained in the closure of its associated global arm  because $b$ is the limit of characteristic points. Hence $\tau(b)\not=\upper(\widetilde{\nu})$ and $b$ is a preperiodic or periodic branch point. In the latter case, $b$ is characteristic and evil.

\smallskip
     $(1.1.)$ We first consider the case that $b$ is preperiodic. Let $\mu=\tau(b)$ be the preperiodic kneading sequence generated by $b$. 
        We claim that $[\overline{\star},\nu]\cap[\overline{\star},\tilde{\nu}]=[\overline{\star},\mu]$.
        
                Let $1\rightarrow S_1 \rightarrow \cdots \rightarrow S_k \rightarrow \cdots$ be the (infinite) internal address of
        $\mu$ and $\mu^k$ be the unique $\star$-periodic sequence associated to
        $1\rightarrow S_1 \rightarrow \cdots \rightarrow S_k$ for $k\in\N$.
        We claim that $b$ is the limit point of characteristic points $p_{S_k}\in T$ with itinerary
        $\tau(p_{S_k})=\upper(\mu^k)$:
        by definition $b$ is the limit point of a sequence of characteristic points $z_{k}$ of period $n_k$.
        We can assume that for all $k\geq l$, the itineraries of $z_k$
        have the same first $n_l$ entries as $\mu$. Fix an entry $S_k$ of the internal address
        of $\mu$.
        Then there are $n_{k_1},n_{k_2}\in\N$ such that $n_{k_1}< S_k \leq n_{k_2}$.
        Since the itinerary $\tau_{k_2}$ of $z_{k_2}$ coincides
        with $\mu$ for the first $n_{k_2}$ entries, $S_k$ is contained in the internal address of $z_{k_2}$.
        Then by~\cite[Proposition 6.6]{BS}, there is a characteristic point $p^\prime_{S_k}$
        with internal address $1\rightarrow S_1 \rightarrow \cdots \rightarrow S_k$ in the
        Hubbard tree of $\upper^{-1}_{n_{k_2}}(\tau(z_{k_2}))$. This forces a characteristic
        point $p_{S_k}\in T$ with $\tau(p_{S_k})=\upper(\mu^k)$.
        By induction on $k$, we get a sequence of characteristic 
         points $p_{S_k}\in T$, which have to converge to
        $b$ because $p_{S_k} \in [z_{k_1}, z_{k_2}]$ for all $k$.
        Since for all $k\in \N$, $T$ contains a characteristic
        point with itinerary $\upper(\mu^k)$ and since these points are contained in the arc
        $[c_0,b[$, we have $\nu > \mu^k$, $\tilde{\nu} > \mu^k$ for all $k\in\N$ and thus, $\nu> \mu$ and $\tilde{\nu} > \mu$.
        
        Assume that $\mu'$ is $\star$-periodic with $\nu, \tilde{\nu} > \mu'$ and let $q' \not=b$ be the characteristic point with
        itinerary $\upper(\mu')$. Since $b$ is the limit of characteristic points,
        Lemma~\ref{LEMcharPointAfterBranchPointNotForced} implies that
        $q'\in [c_0,b[$. Since $p_{S_k}\to b$, there is a $k_0$ such that
        $q'\in \,]c_0,p_{S_{k_0}}[$. By Orbit Forcing~\ref{LEMorbitForcing}, the Hubbard tree of $\mu^{k_0}$
        contains a characteristic
        point with itinerary $\upper(\mu')$ and hence $\mu'<\mu^{k_0}<\mu$. 

\smallskip
    $(1.2.)$ Now suppose that $b$ is periodic.
        Let $\mu$ be the primitive $\star$-periodic kneading sequence such that
        $\tau(b)=\overline{\upper}(\mu)$. Then $\mu\not<\nu$ but $\mu'<\nu$ for all $\mu'<\mu$ by Lemma~\ref{CORevilBranchPointNoUpper}.
        We know that $b$ is the limit point of characteristic points $z_n \in P$. We may assume that
        $c_0<z_n<z_{n+1}<b$.
        Hence, there is a sequence of corresponding characteristic points $\tilde{z}_n \in \widetilde{T}$.
        By compactness of the Hubbard tree and since $\tilde{c}_0<\tilde{z}_n<\tilde{z}_{n+1}<\tilde{c}_1$, this sequence has
        to converge to a point $\tilde{b} \in [\tilde{c}_0,\tilde{c}_1]\subset \widetilde{T}$.
        We know that $\tau(\tilde{z}_n)\to \overline{\upper}(\mu)$. Hence by Lemma~\ref{LEMlimitOfItineraries},
        the itinerary of $\tilde{b}$ is either $\overline{\upper}(\mu)$ or $\mu$, the latter if and only if
        $\tilde{b}=\tilde{c}_1$. In this case, $\mu=\tilde{\nu}$.
     
        In the first case, $\tilde{b}$ is a characteristic point: since $\tilde{b}$ is the limit of characteristic points, all iterates $\iter{f}{i}(\tilde{b})$ are contained in the closure of the global arm of $\tilde{b}$ containing $\tilde{c}_0$. Let $m$ be the period of $\mu$. If $\tilde{b}$ was not $m$-periodic, then there is a $k\in\N$ such that $\tilde{z}_k\in\,]\iter{f}{m}(\tilde{b}),\tilde{b}[\,\subset T_{\overline{\upper}(\mu)}$ and thus $\tau(\tilde{z}_k)=\tau(\tilde{b})=\overline{\upper}(\mu)$, which is not true.
        We claim that either $\mu\leq\tilde{\nu}$ or $\mu\not<\tilde{\nu}$ and $\mu'<\nu$ for all $\mu'<\mu$.
        If $\tilde{b}$ is an inner point, then Theorem~\ref{THMexistenceBifuracation}
        guarantees the existence of a characteristic point
        with itinerary $\upper(\mu)$ or $\mu$ and hence $\mu\leq\tilde{\nu}$. If $\tilde{b}$ is an evil
        branch point, then again by Lemma~\ref{CORevilBranchPointNoUpper} we have
        that $\mu\not<\tilde{\nu}$ and that $\mu'<\tilde{\nu}$ for all $\mu'<\mu$.
        Observe that there is no $\mu'\not<\mu$ such that both $\mu'<\nu$ and
        $\mu'<\tilde{\nu}$:
        for any $\mu'<\nu$ and $\mu'\not<\mu$ there is a characteristic periodic
        point $z\in\,]b,c_1]\,\subset T$ with itinerary $\upper(\mu')$.
        Since $b$ is characteristic itself, Lemma~\ref{LEMcharPointAfterBranchPointNotForced}
        says that there is no characteristic point in $\widetilde{T}$ with itinerary $\tau(z)$, and so
        $\mu'\not<\tilde{\nu}$.

    \item The second possibility is that $a \in \,]c_0,b[$. We first show that $a\in P$:
       the arc $]a,b[$ contains a precritical point because $a$ and $b$ have different itineraries:
        by minimality of $T$, $b$ cannot have the same itinerary as any point of $P$. Let $I$ be the set of points in $[c_0,b]$ with itinerary $\tau(b)$. Then $I=\,]\xi,b]$ with $\xi$ precritical by Lemma~\ref{LEMptsOfSameItinerary}. So, if $a\in I$,
        then there is an open neighborhood $U$ of $a$ such that $U\subset I$ and hence $P\cap I\not=\emptyset$.
        If $a$ is a trivial limit, then clearly $a\in P$. Otherwise, let $\xi$ be the precritical point in $]a,b[$
        such that there is no precritical point in $]a,b[$ with lower step and set $k:=\,${\sc step}$(\xi)$.
        Since $a$ is a limit point of characteristic points, we have that for all $l\in\N$, 
        $\iter{f}{l}(a)\not\in \,]a,c_1[$. Hence, $\iter{f}{k}([a,\xi])$ covers $[a,\xi]$ homeomorphically.
        Therefore, there is a $z\in [a,\xi[$ with $z=\iter{f}{k}(z)$.
        If $z=a$, we are done.
        Otherwise, let $w\in [a,\xi]$ be the periodic point with lowest period.
        If $w\not=a$ (otherwise we are done), then the point $w$ is characteristic:
        if not, then there is an $s<k$ such that $\iter{f}{s}([a,w])$ covers
        $[a,w]$ homeomorphically. But this yields a periodic point
        $\tilde{w}\in\,]a,w[$ of period smaller than the one of $w$, a contradiction to the choice of $w$.
        Hence $w\in P$ and $a\not=\sup(P)$.
        This proves that $a$ is periodic and as a limit point of characteristic points it is characteristic
        itself.

        Let $\tau$ be the itinerary of $a$, $n$ be its period and let $\mu$ be the $\star$-periodic kneading sequence
        such that  $\tau=\upper(\mu)$ or $\tau=\overline{\upper}(\mu)$ according as $a$ is tame or not. Let $q$ be the number of arms at $a$ and $Q$ the period of the local arm at $a$ pointing to $c_1$.
         By Lemma~\ref{LEMcharPointsAreForced}, there is a characteristic point $\tilde{a}\in\widetilde{T}$ which has itinerary $\tau$ and which is of the same type as $a$. Let $\tilde{q}$ be the number of arms at $\tilde{a}$ and $\tilde{Q}$ the period of the local arm at $\tilde{a}$ pointing to $\tilde{c}_1$. Then $q=\tilde{q}$ if and only if $Q=\tilde{Q}$.
        Theorem~\ref{THMexistenceBifuracation} implies
        that $T$ contains a characteristic point $x$ with itinerary $B_Q(\tau)$ (it is not possible that $x=c_1$), and $\widetilde{T}$
        contains a characteristic point $\tilde{x}$ that has itinerary $B_{\tilde{Q}}(\tau)$ or $\upper^{-1}_{\tilde{Q}n}(\tau)=B_{\tilde{Q}}(\mu)$.
        In the last paragraph of the proof of Lemma \ref{LEMcharPointsAreForced}, we have seen that a Hubbard tree cannot contain two characteristic points with itinerary $B_Q(\tau)$ and $B_{\tilde{Q}}(\tau)$  respectively for $Q\not=\tilde{Q}$. This together with Lemma~\ref{LEMcharPointAfterBranchPointNotForced} yields that if $q\not=\tilde{q}$
and $\mu'$ $\star$-periodic with $\mu'<\nu,\,\mu'<\tilde{\nu}$, then $\mu'\leq\mu$.
If $\tau=\upper(\mu)$, then $\mu<\nu$ and $\mu<\tilde{\nu}$. If $\tau=\overline{\upper}(\mu)$, then $\mu$ is primitive and
depending on whether $a$ ($\tilde{a}$ respectively) is an inner or branch point, $\mu<\nu$ ($\mu<\tilde{\nu}$) or $\mu\not<\nu$ and  $\mu'<\nu$ for all $\mu'<\mu$ ($\mu\not<\tilde{\nu}$ and $\mu'<\tilde{\nu}$ for all $\mu'<\mu$). In all cases, $\mu$ is the wanted kneading sequence of the Branch Theorem. 

              If $q=\tilde{q}$, then $[\overline{\star},\nu]\cap[\overline{\star},\tilde{\nu}]=[\overline{\star},B_Q(\mu)]$: let $\mu'$ be $\star$-periodic such that $\mu'<\nu$ and $\mu'\not\leq B_Q(\mu)$. There is a characteristic point $y\in T$ with $\tau(y)=\upper(\mu')$ and $y\in\,]x,c_1[$. Since $p\not=x$, Lemma~\ref{LEMcharPointAfterBranchPointNotForced} yields that there is no characteristic point in $\widetilde{T}$ with itinerary $\tau(y)$ and thus, $\mu'\not<\tilde{\nu}$.           
\qedhere   \end{enumerate}
\end{proof}

Using the notion of wakes, we can interpret the Branch Theorem the following way:
      let $\nu\not=\tilde{\nu}$ be two $\star$- or preperiodic kneading sequences.
        Then either one kneading sequence is contained in a (non-admissible or admissible) subwake
        of the other, or there is a third $\star$-periodic or preperiodic kneading sequence $\mu$
        such that $\nu$ and $\tilde{\nu}$ are contained in two different (non-admissible or
        admissible) subwakes of $\mu$. If one of the two given kneading sequences is contained in a
        non-admissible subwake of $\mu$ (or of $\nu$, $\tilde{\nu}$), then $\mu$ (or $\nu$, $\tilde{\nu}$) is $\star$-periodic and primitive.

\medskip
    Suppose the Hubbard trees of $\nu$ and $\tilde{\nu}$ are non-admissible and they contain
    an evil periodic branch point $b$ and $\tilde{b}$ such that $\tau(b)=\tau(\tilde{b})$.
    Then $\nu$ and $\tilde{\nu}$  branch off into a non-admissible wake at the same admissible kneading sequence $\mu$:
    an immediate consequence of the Orbit Forcing~\ref{LEMorbitForcing} is that the
    evil characteristic points closest to $c_0$ in $T$ and $\widetilde{T}$ have the same
    itinerary. The $\star$-periodic kneading sequence $\mu$ associated to these evil characteristic
    points is admissible such that $\mu\not<\nu$ and $\mu\not<\tilde{\nu}$ and for all
    $\tilde{\mu}<\mu$ we have both $\tilde{\mu}<\nu$ and $\tilde{\mu}<\tilde{\nu}$ by~\cite[Proposition 6.5]{BS}.


\begin{thebibliography}{MM}
    \bibitem[BS]{BS} H. Bruin, D. Schleicher, \emph{Symbolic dynamics of quadratic
                polynomials}, Preprint, Institut Mittag-Leffler (2002).


    \bibitem[DH]{orsay} A. Douady, J.H. Hubbard, \emph{Etude dynamique des polyn\^{o}mes
                complexes}, Publication math\'{e}matiques d'Orsay 84-02 {\bf 75} (1984) and 85-04 {\bf 138} (1984).

   
    \bibitem[K]{kauko2} V. Kauko, \emph{Shadow trees of Mandelbrot sets},
                Fundamenta Mathematicae {\bf 180}, 35--87 (2003).

    \bibitem[Ka]{kafflThesis} A. Kaffl, \emph{Hubbard trees and kneading sequences for unicritical and cubic polynomials}, Ph.D. Thesis, International University Bremen (2006).
    
    \bibitem[L]{lavaurs} P. Lavaurs, \emph{Une d\'{e}scription combinatoire de l'involution definie par M sur les rationnels
                \`a d\'enominateur impair}, Comptes Rendus Acad. Sci. Paris {\bf 303}, 143--146
                (1986).

    \bibitem[LS]{lausch} E. Lau, D. Schleicher, \emph{Internal addresses in the Mandelbrot set
                and irreducibility of polynomials}, Preprint,
                Institute for Mathematical Sciences, Stony Brook (1994).

    \bibitem[N]{nadler} S. Nadler, Jr, \emph{Continuum theory --- an introduction}, Marcel Dekker (1992).

    \bibitem[M] {milnor} J. Milnor, \emph{Periodic orbits, external rays and the Mandelbrot set;
                    an expository account}, Ast\'{e}risque {\bf 261}, 277--333 (2000).

    \bibitem[Pe]{penrose} C. Penrose,  \emph{On quotients of shifts associated with dendrite Julia sets of
                quadratic polynomials}, Ph.D. Thesis, University of Coventry (1994).
                
    \bibitem[Po]{poirier} A. Poirier, \emph{On post critically finite polynomials; part two: Hubbard trees},
                Preprint, Institute of Mathematical Sciences, Stony Brook (1993).

    \bibitem[S]{sFibers} D. Schleicher, \emph{On fibers and local connectivity of Mandelbrot and Multibrot
                sets}, In: M. Lapidus, M. van Frankenhuysen (eds): \emph{Fractal Geometry and Applications: A Jubilee
                of Benoit Mandelbrot.} Proceedings of Symposia in Pure Mathematics (AMS), to appear.

    \bibitem[T]{thurston} W. Thurston, \emph{On the geometry and dynamics of iterated rational maps},
                Preprint, Princeton University (1985).

\end{thebibliography}
\end{document}